\documentclass[letterpaper,11pt]{amsart}
\usepackage[margin=1.25in]{geometry}
\usepackage[maxbibnames=99,maxcitenames=99]{biblatex}
\usepackage{graphicx}

\usepackage{amssymb}
\usepackage{amsmath}
\usepackage{mathtools}

\newcommand{\Z}{{\mathbb Z}}
\newcommand{\N}{{\mathbb  N}}

\newcommand{\lra}{\leftrightarrow}
\newcommand{\dlra}{\leftrightarrow}

\def\PP {{\mathbb P}}

\newtheorem{Theorem}{Theorem}

\newtheorem {lemma} [Theorem]    {Lemma}

\newtheorem {proposition}[Theorem]    {Proposition}

\newtheorem{remark}[Theorem]{Remark}

\newcommand {\eps}{\varepsilon}

\begin{document}

\title{On loops in critical high-dimensional percolation}
\author{Amelia Carpenter \and Wendelin Werner}
\address{DPMMS, University of Cambridge}

\begin {abstract}
We show the following results about critical Bernoulli percolation in high dimensions:
In the box of side-length $N$, there exist self-avoiding open loops of diameter comparable to $N$, and the collection of these self-avoiding loops has a non-trivial scaling limit (if viewed in the Hausdorff topology) as $N$ tends to infinity.
This feature contrasts with the  proliferation of ``typical'' percolation clusters pointed out by Michael Aizenman almost three decades ago. In other words, we show that among the many large clusters in a given box, only a handful will contain a self-avoiding loop of diameter greater than a fixed fraction of the side-length of the box.
\end {abstract}

\maketitle

\bigbreak

\null
{\centerline {\em
This article is dedicated to Michael Aizenman on the occasion of his 80th birthday.}
}
\vskip 1cm

\section {Introduction}

\subsection {General introduction}

One way to describe the phase transition for bond Bernoulli percolation in $\Z^d$ goes as follows: Consider the box $\Lambda_N := [-N,N]^d$ viewed as a subgraph of $\Z^d$, erase or keep each of its edges independently with probability $1-p$ and $p$ respectively, consider the family $(C_j)_{j \in J}$ of the connected components of the resulting subgraph of $\Lambda_N$. Its rescaled version $(K_j:= C_j / N)_{j \in J}$ is then a random family ${\mathcal K}= {\mathcal K}(N,p,d)$ of disjoint compact subsets of $[-1, 1]^d$. If one fixes $p$ and looks at the behaviour of this family when $N \to \infty$, then for any $d \ge 2$ (there are many textbooks presenting this, see for instance \cite {G}), there exists a value $p_c = p_c (d) \in (0,1)$ such that:  (a) when $p< p_c$, all the clusters of ${\mathcal K}$ will be very small (with a probability that goes to $1$ as $N \to \infty$), while (b) when $p>p_c$, there will be one large ``dense'' cluster and otherwise very small ones (with a probability that goes to $1$ as $N \to \infty$). This raises of course the question about what happens when $p$ is equal to this critical value $p_c$ [which is the model called {\em critical percolation}]. In that case, the behaviour of ${\mathcal K}$ as $N \to \infty$ will depend on the dimension. For small values of $d$, it is believed that the law of ${\mathcal K}$ converges to that of a random collection of compact sets in $[-1,1]^d$ with the property that for any small $a$, the number of clusters of diameter at least $a$ converges to that of some non-trivial finite random variable. This is usually referred to as the existence of a {\em scaling limit} -- and a number of results in this direction have been obtained over the years (but many of the fundamental statements are still conjectures).
On the other hand, when $d$ is large (typically $d \ge 7$ -- this is sometimes referred to as $d$ being above the critical dimension) which is the  focus of the present paper, this is no longer believed to be true.

The study of critical percolation-type lattice models above their critical dimension has a long distinguished history including the two-point function estimates by Hara and Slade \cite {HS} via the lace expansion technique (that we will recall in the next section -- see the monographs \cite {Sladebook,HH} on the topic for an extended list of references including \cite {H,S,FH,FH2}) and it has also been the topic of a number of recent and ongoing works including \cite {HHHM, CH, CHS, DCP, CFHJ, CCHS0, CCHS, CCHS2, ASS, HMS, }.
One main line of results and conjectures in this case is that the number of large clusters (i.e. of diameter comparable to the size of the box $\Lambda_N$) will tend to infinity (when $N \to \infty$). More specifically, in the previous setup, the clusters in ${\mathcal K}$ with diameter at least $a$ (when $a$ is fixed and small, while $N$ grows) do proliferate -- their number would in fact typically be of the order $N^{d-6 + o(1)}$. Furthermore, these clusters would be tree-like (i.e., the self-avoiding loops that they contain are all of size much smaller than $N$) and would in fact in some sense resemble the trace of integrated super-Brownian excursions. In particular, a large cluster $C_j= N K_j$ in $[-N, N]^d$ would have circa $N^4$ points and the paths joining two far-away points within such a $K_j$ would look like Brownian paths. This loose description can be made more precise and in many cases proven --  this is what some of the aforementioned activity in this field has  been about. The key to deriving these results is to first obtain estimates about the ``two-point function'', i.e. the probability that two far-away points belong to the same cluster. With such estimates in hand, one can use a combination of diagrammatic moment bounds, the BK and FKG inequalities, the Paley-Zygmund inequality to deduce the aforementioned results for Bernoulli percolation. The proliferation of these large clusters (building on two-point function estimates) has been first pointed out and investigated by Aizenman \cite {Ai} in the mid-1990s.

\begin{figure}[h]
  \centering
  \includegraphics[width=.6\textwidth]{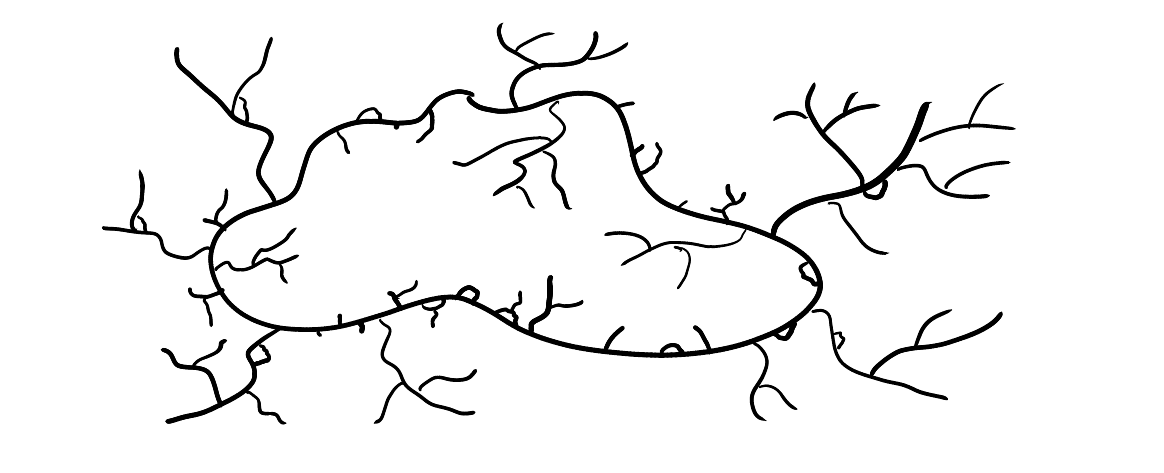}
  \caption{\label{p0}Impression of a cluster containing a large loop}
  \end{figure}

The main purpose of this paper is to show that, while the previous description in terms of trees is that of {\em typical} and proliferating large clusters, there will exist {\em exceptional} large clusters (at any scale) that are not macroscopically tree-like (i.e., that do contain large self-avoiding loops as sketched in Figure~\ref {p0}). In fact, for any fixed $a$, the number of clusters that contain loops of diameter greater than $aN$ in $\Lambda_N$  will be tight, and the largest diameter of an open self-avoiding loop contained in $\Lambda_N$ will be of order $N$. In other words, among the aforementioned $N^{d-6}$ clusters, a tight number will actually contain loops of diameter greater than $aN$ for small fixed $a$. Also, each of these clusters will (in some sense that we will explain) in fact essentially contain only one big loop (i.e., any two big self-avoiding loops in these clusters will be $o(N)$-close in Hausdorff distance). We will also show the existence of (subsequential) scaling limits for these collections of large loops. So, there are scaling limits after all if one focuses on clusters that are topologically non-trivial!

As we shall discuss, these loops should asymptotically look like Brownian loops (and the corresponding cluster will otherwise contain trees attached to this loop) --  the collection of all these large loops then would give rise to Brownian loop-soups (as defined in \cite {LW} -- this is a Poisson point process of Brownian loops). This is similar to the relation between the backbone of the infinite incipient cluster in high dimensions and Brownian motion (see \cite {HS2,HS3,HHHM} among other papers).

While this type of result is new and possibly at first glance surprising, the techniques at work to derive tightness, existence and properties of such large loops rely to a large extent on the tool-box developed in the many papers on the topic of high-dimensional percolation i.e., one uses mainly the same type of diagrammatic bounds, BK inequalities and Paley-Zygmund second moment ideas that have been developed an used in many of the papers on the subject.
\medbreak

\begin{figure}[h]
  \centering
  \includegraphics[width=\textwidth]{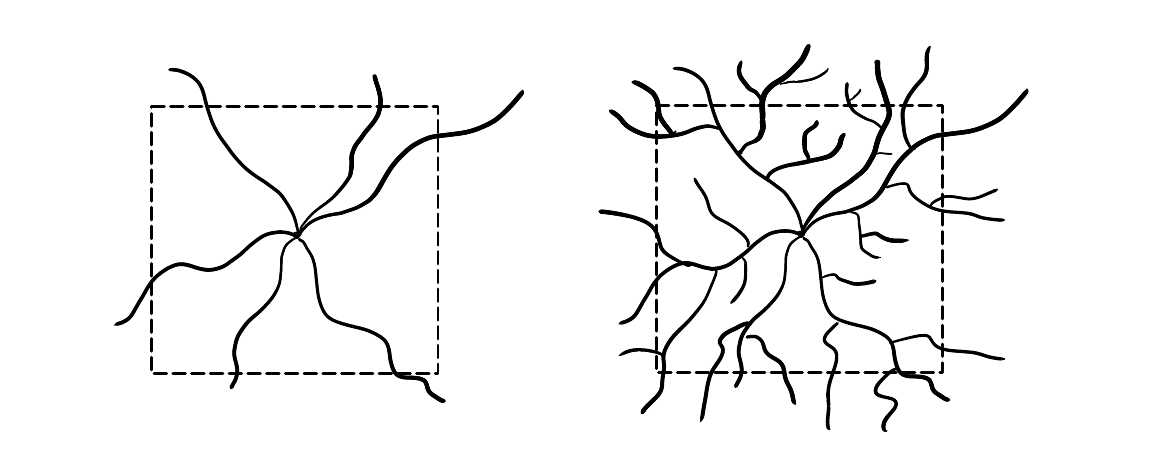}
  \caption{\label{p1}A $6$-arm starfish and a sketch of a cluster containing such a starfish. These clusters are tight when the dimension of space is $12$}
  \end{figure}

We will also explain that other special types of exceptional clusters (or sets of touching clusters) will give rise to tight families in the large scale limit. A particular natural one is depicted in Figure \ref {p1} -- these are the clusters that contain starfishes with $l$ arms (i.e., there exist $l$ disjoint open paths of diameter at least $aN$ in the cluster that all start at the same point), that turn out to be tight in the box $\Lambda_N$ when the spatial dimension $d$ is even and equal to $2l$ (since we anyway assume that $d \ge 7$, we need that $l \ge 4$). This will be Proposition~\ref {starfish} in Section \ref {S73}.

Let us mention that our results seem to be quite a different story from the interesting case of percolation on high-dimension tori, that has been the subject of a number of papers (including recent ones -- see for instance \cite {HStori,HMS,BN} and the references therein) and that may at first glance look similar. The situation in the torus is indeed different as the macroscopic cluster and its non-contractible cycles will be created by the connection between what would correspond to ``typical very very big clusters'' in a much larger box for percolation in the universal cover (so that the phenomenon at work in tori is more related to the emergence of the giant component in the complete graph).

Finally, let us stress that like many of the papers on scaling limits in these last decades, {\bf this paper is in more than one way ``in the spirit'' of ideas put forward by Michael Aizenman to whom this paper is dedicated on the occasion of his 80th birthday}: He advocated the study of scaling limits (albeit below the critical dimension), the use of tightness ideas such as in \cite {AB}, and he also showed the way (then exploited by others, including of course his own PhD students who form a substantial fraction of the bibliography) on how to exploit diagrammatic bounds and the triangle conditions \cite {AN,Ai} in order to study high-dimensional critical and near-critical percolation, and he is the one who started the study of geometric features and proliferation of critical percolation clusters in high dimensions in \cite {Ai}.

\subsection {Precise setup and main statements}
Let us now state the three main results about loops in high-dimensional percolation that we will derive.
We consider critical bond percolation in $\Z^d$ for $d$ sufficiently large -- so each edge is open or closed independently with probability $p_c$, where $p_c$ is the critical threshold described above (which is also the one above which infinite open clusters appear; see for instance \cite {G} for basics about the model, or the BK inequality that we will repeatedly use).

{\bf Throughout the paper, we are going to work under the assumption that $d \ge 7$ and that the following bound on two-point connectivity probabilities holds as $|x-y| \to \infty$:}
\begin {equation}\label{tpe}
P [ x \leftrightarrow y ] \asymp \frac 1 {|x-y|^{d-2}}
\end {equation}
(here and in the sequel, $x \lra y$ stands for the event that there is a path of open edges connecting $x$ and $y$, and  $\asymp$ will stand for the fact that the ratio between the left-hand side and the right-hand side is bounded and bounded away from $0$). For nearest-neighbour Bernoulli percolation, this is believed to be true for all $d \ge  7$, and known to hold for all sufficiently large $d$ (more precisely, this has been shown via the Brydges-Spencer lace expansion techniques when $d \ge 19$ by Hara and Slade in \cite {HS} and extended to $d \ge 11$ by Fitzner and van der Hofstad \cite {FH} with computer assisted computations), see the monographs \cite {Sladebook,HH} for an overview. So our results are unconditional results when $d \ge 11$, and for $d \in [7,10]$, they are conditional on the conjecture that \ref {tpe} holds in these dimensions.  For spread-out percolation, this is known to hold for all $d \ge 7$ (see \cite {HHS}), and our proofs could be adapted to work in this case, but for rather obvious presentation purposes, we will stick to the nearest-neighbour Bernoulli percolation in this paper.
A new powerful and different (i.e., without the lace expansion technique) approach to the derivation of estimates for two-point functions in high-dimensional sufficiently spread-out percolation has been pointed out in \cite {DCP}. Finally, to conclude this paragraph, let us note that in fact (see \cite {H,FH,FH2} and the references therein) the lace expansion technique can in fact  be pushed to show that $P [ x \lra y ] \sim c_d / |x-y|^{d-2}$ for some constant $c_d$ when $d \ge 11$, but will not rely on this stronger result here.

For $x \in \Z^d$, we define $L(x)$ to be the set of sites $y$ in $\Z^d$ such that
$$x \Leftrightarrow  y := ( x \leftrightarrow y ) \circ ( x \leftrightarrow y )$$
holds (here, we use $A \circ B$ to denote the event that $A$ and $B$ are realized disjointly, as customary in the statements of the BK inequality). In other words,  $L(x)$ is the set of points such that there exists a loop that does not use the same edge twice (mind that the loop is allowed to go twice through a given site, so a figure eight type loop is for instance allowed) that goes through both $x$ and $y$. By abuse of terminology, we will call such loops to be self-avoiding
throughout this paper (they are ``edge''-self-avoiding), and we will refer to the loops that do not visit any site twice as ``site-self-avoiding''.

\begin{figure}[h]
  \centering
  \includegraphics[width=
  \textwidth]{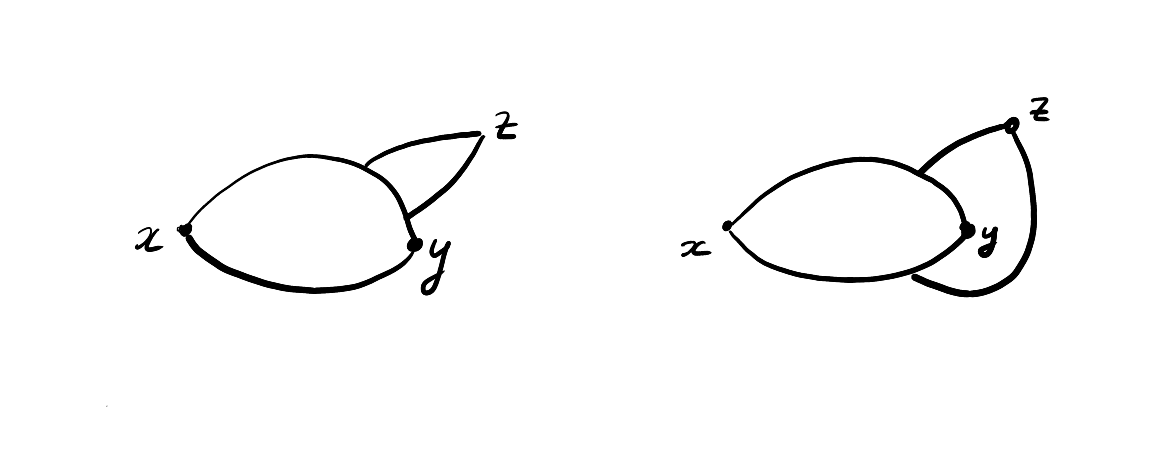}
  \caption{\label{f1}The relation $\Leftrightarrow$ is transitive: Depiction of the two possibilities when $x \Leftrightarrow y$ and $y \Leftrightarrow z$ both hold -- in either case, one has $x \Leftrightarrow z$}
  \end{figure}

It is straightforward to check that if $x \Leftrightarrow  y$ and $y \Leftrightarrow  z$, then one can find two disjoint open self-avoiding paths that join $x$ and $z$: One can first draw the two self-avoiding disjoint paths $l_1$ and $l_2$ from $x$ to $y$, and then consider the parts $l_1'$ and $l_2'$ of the two disjoint open self-avoiding paths from $z$ to $y$ up to their respective first hitting of $l_1 \cup l_2$, and then choose the remainder of the connections from $z$ to $y$ by choosing the appropriate parts of $l_1 \cup l_2$ depending on where $l_1'$ and $l_2'$ do hit $l_1 \cup l_2$ (see Figure \ref {f1}). It therefore follows that $\Leftrightarrow $ is an equivalence relation.

\medbreak
{\sc Definition.}
{\em Equivalence classes for $\Leftrightarrow$ are called loop-clusters.}
\medbreak
So, a loop-cluster is a maximal connected union of open self-avoiding loops, and  $L(x)$ is the loop-cluster that contains $x$.
Making sense and defining these loop-clusters is somewhat similar to finding a workable definition for {\em the backbone} (starting from a given point $x$) of the ''incipient infinite cluster'' (the analog would be that it is the set of points $y$ for which $( x \leftrightarrow y ) \circ ( y \leftrightarrow \infty )$ for this modified percolation model).
\medbreak
Note that the min-flow max cut property gives the following equivalent approach: $x \Leftrightarrow y$ if and only if the removal of any single open edge does not disconnect $x$ from $y$ (and the transitivity then becomes immediate and does not require drawing Figure \ref {f1}). In Graph Theory, our loop-clusters are  referred to as 2-connected components. As we shall see, in our percolation context, the definition via max flow rather than min cut appears to be the most natural one (it allows the use of the BK inequality for instance).

\begin{figure}[h]
  \centering
  \includegraphics[width=.8\textwidth]{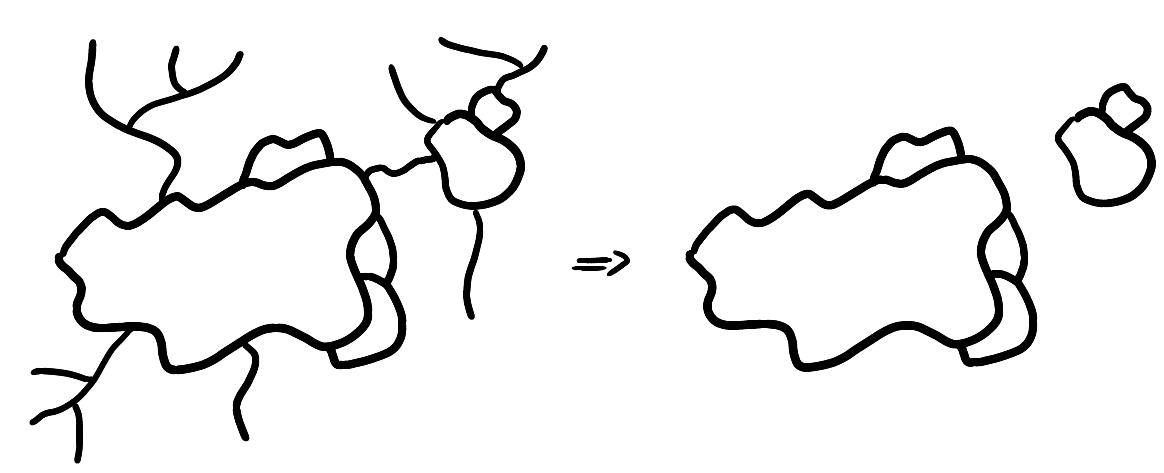}
  \caption{\label{w1}A percolation cluster and its two loop-clusters}
  \end{figure}

Note of course that each loop-cluster is contained in a percolation cluster but that a percolation cluster can contain several disjoint loop-clusters (see Figure \ref {w1}).
Note also that any loop-cluster contains a self-avoiding loop of the same diameter (because any two points in the loop-cluster will be part of a self-avoiding loop that goes through both, and this self-avoiding loop is then part of the loop-cluster).

If we restrict the percolation to a subdomain $\Lambda$ of $\Z^d$, we can similarly define the loop-clusters within $\Lambda$, and define $L_\Lambda (x)$ to be the loop-cluster within $\Lambda$ that contains $x$. In most of the paper, we will be considering the case where $\Lambda = \Lambda_N = [-N,N]^d$.
\medbreak

The next sections will be mostly devoted to explain why the following two propositions hold; natural variants do hold as well -- one can for instance easily replace the collection of loop-cluster in $\Lambda_N$ by the collection of loop-clusters of $Z^d$ that do intersect $\Lambda_N$.

\begin {proposition}[Tightness and existence of large loop-clusters]
\label{mainprop}
\label {main1}
 For any $a<1$, the number $n(a,N)$ of loop-clusters in $\Lambda_N$ of diameter greater than $aN$ is tight (as $N \to \infty$). Furthermore, for any $k$ and $\eta$, one can find $a$ small enough such that $P [ n(a,N) \ge k ] \ge 1- \eta$ for all large enough $N$.
\end {proposition}

This indicates that large loop-clusters exist but do not proliferate.
More precisely (for each given $N$), let us define $(R_k^N)_{k \ge 1}$ to be the diameters of the loop-clusters in $\Lambda_N$ listed in decreasing (i.e. non-increasing) order. Then, as a collection of random sequences in a given compact interval, the family of laws of $(R_k^N / N)_{k \ge 1}$ is always tight (in the topology given by the convergence of any finite-dimensional marginals) when $N$ varies. One way to reformulate  the proposition is that any subsequential limit is supported on the set of sequences $(r_k)_{k \ge 1}$ of strictly positive numbers with $\lim_{k \to \infty} r_k = 0$ (strictly positive corresponds to the existence, while the fact that $r_k \to 0$ corresponds to the tightness part of the proposition).

\begin{figure}[h]
  \centering
  \includegraphics[width=.8\textwidth]{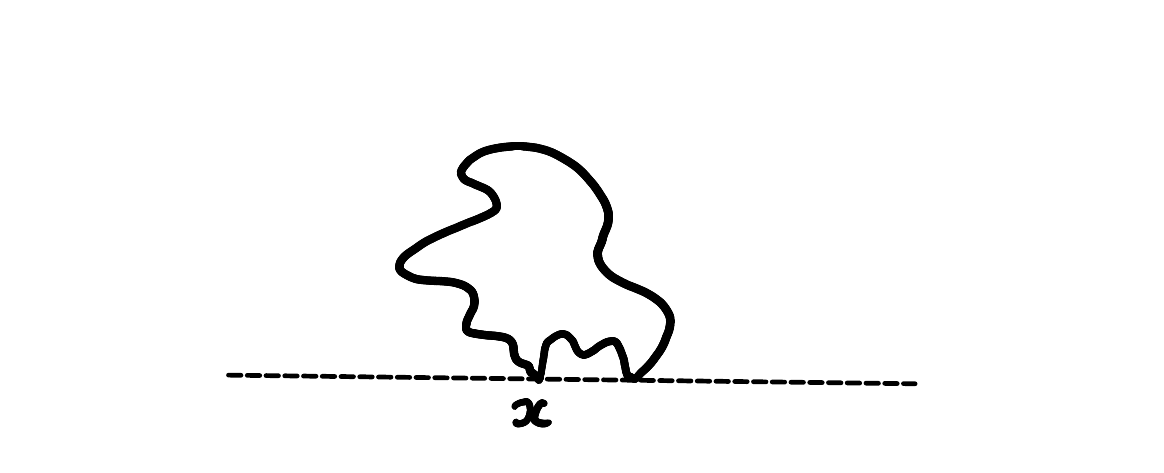}
  \caption{\label{fig:6a} ``Lowest'' point of a large loop}
  \end{figure}

Let us outline how one possible proof of such a tightness result goes (a more high-level ``reason'' for the result to hold is the loop-soup result from \cite {LupuW} that we will recall at the end of this section); this also indicates how it relates to two-point functions in half-spaces that have been the focus of a number of recent papers that we will mention in the present paper:  If an open self-avoiding loop in $\Lambda_N$ has $L^\infty$-diameter greater than $a  N$, then we can look at its point (or one of its points if there are several) with minimal first coordinate. From such a point $x$, the loop will be a loop of large diameter in the half-space consisting of the points with first coordinate larger than that of $x$.

Now, the boundary to boundary two-point functions in half-spaces by Chatterjee-Hanson \cite {CH} suggest that for each given $x$, the probability of such an open self-avoiding loop existing will be of order $N^{-d}$.
So, when summing over the $O(N^d)$ points $x$, we get that the expected number of such extremal points (for fixed $a$, and all $N$) is indeed bounded and bounded away from $0$. This upper bound already shows that the number of big loop-clusters is tight (since there will be anyway at least one such point per big loop-cluster). To complete the argument (i.e., to see  that such large loop-clusters do indeed exist when $a$ is small enough) which loosely amounts to checking that large loop-clusters will tend to have very few such minimal point of self-avoiding loops, the (classical) idea is to show that the second moment of the number of minimal points is bounded as well (which by Cauchy-Schwarz will imply that the probability of there existing at least one such minimal point is bounded away from $0$ independently of $N$).

\medbreak

We move to our second main result, which loosely speaking says that large loop-clusters have essentially only one large self-avoiding loop -- this is the idea conveyed by the sketch in Figure \ref {p0}. In this statement, $\alpha$ is any constant in $( 2/(d-4), 1)$ (and the smaller $\alpha$ is, the stronger the statement becomes):
\begin {proposition}[One ``large loop'' per large loop-cluster]
\label{main2}
For all fixed $a < 1$,  the  probabilities of the following two events for critical percolation in $\Lambda_N$ are going to $1$ as $N \to \infty$:
\begin {itemize}
 \item
Any two open self-avoiding loops $\gamma_1$ and $\gamma_2$ with diameter greater than $aN$ that are part of a same percolation cluster do intersect and the Hausdorff distance between $\gamma_1$ and $\gamma_1 \cap \gamma_2$ (and also between $\gamma_1$ and $\gamma_2$) is bounded by $N^\alpha$.
\item
Any open (edge)-self-avoiding loop $\gamma$ with diameter greater than $aN$ contains a ``site-self-avoiding'' loop $\tilde \gamma$ that is at Hausdorff distance smaller than $N^\alpha$ from $\gamma$.
\end {itemize}
\end {proposition}

There are alternative or stronger statements in the same vein (for instance in terms of presence of ``pivotal'' edges along each loop).
 We will also explain why the typical number of points in such large loop-clusters described in Proposition \ref {mainprop} will be of order $N^2$. This, as well as Proposition \ref {main2} itself, is  reminiscent of (and closely related to) results about the structure of the backbone of the incipient infinite cluster.

 We will also show (this is going to be Lemma \ref {L:no-pinching}) that the large self-avoiding loops will not be ``close to pinching'' which alongside some analytical topology
 considerations and Propositions~\ref {mainprop} and \ref {main2} will imply the existence of subsequential limits for rescaled large loop-clusters (in the weak topology associated to the Hausforff topology), where the limit is supported on the set of continuous self-avoiding loops in $[-1,1]^d$. More specifically, let $(L_j(N), j \ge 1)$ denote the collection of all loop-clusters in $\Lambda_N$, renormalized by a fact $1/N$ so that they are now subsets of $(\Z/N)^d \cap [-1, 1]^d$ and ordered by decreasing diameter (using some deterministic rules to break possible ties). Then:
 \begin {proposition}[Subsequential limits]
 \label {main3}
For any $N_n' \to \infty$, there exists a subsequence $N_n \to \infty$ and a random collection of self-avoiding disjoint continuous loops $(L_j, j \ge 1)$ in $[-1, 1]^d$ with the property that the diameter of $L_j$ almost surely tends to $0$ as $j \to \infty$, so that the law of $(L_j (N), j \ge 1)$ converges to that of $(L_j, j \ge 1)$ [here, the topology is that induced by the Hausdorff distance for finite-dimensional marginals].
 \end {proposition}

We plan to explain in subsequent work why the law of this scaling limit has to be that of  a Brownian loop-soup (as defined in \cite {LW})
in the scaling limit (the Brownian aspect is somewhat similar to the fact that the IIC backbone converges to a Brownian motion such as in \cite {HHHM} for instance, and the Poissonian part follows from the independence properties of percolation). It is not yet entirely clear to us whether such a loop-soup would necessarily be critical (i.e., of the exact intensity that possesses the rewiring property and is closely related to the Gaussian Free Field, see e.g. \cite {WP} and the references therein) -- but there seems to be very strong arguments that suggest it.

One such piece of evidence comes from the corresponding new results about
 percolation of Brownian loops on cable-graphs in high dimensions \cite {LupuW} -- the ideas of that paper actually prompted the investigation that led to the present paper.
 It indeed follows from the switching identities derived in (the July 2025 version of) \cite {W2025} that if one considers a  Brownian loop-soup on the cable-graph of $\Lambda_N$ and keeps only the Brownian loops of size smaller than $N^{1- \eps}$ for some positive fixed small $\eps$, then there will exist large clusters (i.e., clusters of this Poissonian collection of Brownian loops of size much smaller than $N$) that will contain large self-avoiding cycles, and the switching-type identities from \cite {W2025} can then be also used to show (without technical estimates!) that when $d \ge 7$, the collection of these large cycles (rescaled by $N$) will converge (in distribution) to a critical continuum Brownian loop-soup in $[-1,1]^d$  when $N \to \infty$ -- see \cite {LupuW} for a detailed discussion and the proof (this result is essentially the intensity doubling conjecture formulated by Lupu in \cite {Lupu}, building also on some related work by Kassel and L\'evy \cite {KL}). So, the loop-soup with cut-off provides a particular spread-out percolation model for which this convergence of large cycles to the Brownian loop-soup is actually proved.

\medbreak

Another remark is that all our proofs are based on properties of the two-point functions (in the whole space and in the half-space) that are known to hold \cite {HH,CHS,HMS} also in the near-critical regime, so that the results do hold in that setting as well.
In other words, if we fix some positive $\kappa$, and choose a sequence $p(N)$ with $p(N) \in [p_c - \kappa / N^2, p_c + \kappa/ N^2]$, then Propositions \ref {mainprop}, \ref {main2} and \ref {main3} do still hold if for each $N$ one replaces percolation with parameter $p_c$ by percolation with parameter $p(N)$.
We will briefly comment on the loop-soup conjecture in that setting in the final section.

\section {Review of some existing estimates in half-spaces and related facts}
\label {S2}
\label{S21}

Let us now review and recall some results about connection probabilities in half-spaces and the related question of
flow out of boxes. The goal is here just to present only the pieces of the puzzle that we will use rather than the strongest existing statements. Indeed (we will comment on this in Section \ref {Scomment}), many of these upper bounds turn out to be sharp. Most of these result come from the paper \cite {CH} by Chatterjee and Hanson. They are also closely related to the results and considerations of the paper \cite {CHS} by Chatterjee-Hanson-Sosoe and the recently posted paper \cite {panisschapira2} by Panis and Schapira. Of course, there are a large number of papers devoted to the study of aspects of critical Bernoulli percolation in high dimensions and the infinite incipient cluster such as \cite {CC,HJ,KN,KN2,CCHS0}, and many of the considerations that will be discussed here are in some sense related to them as well.

Most of the results of these papers that we will use are theorems and propositions on their own right, but since we will use them here as tools to derive our results, we will refer to them as lemmas (hoping that this does not offend anyone!). Also (and we will comment further on this -- see for instance just after Lemma \ref {L4}), we will only state the parts of these results that we will actually need, in order to clarify what we really build upon.

In this section and throughout this entire paper, we will use $C$ to denote constants that might change from one equation to the other --- when more than one constant are involved in one equation, we will use denote them by $C$, $C'$ etc.

\subsection {The half-space bound}

The first main statement that we want to mention here is the connection probability in the
half-space $H := \N \times \Z^{d-1}$. Here and in the sequel, we will use the standard notation $ x \lra_\Lambda y$ to  denote the event that $x$ and $y$ are connected by a path of open edges that stays in a set $\Lambda$. $0$ will denote the origin in $\Z^d$.

\begin {lemma}[\cite{CH}]
\label {L4}
There exists a constant $C$ such that for all $x \in H$, $P [ 0 \lra_H x ] \le C / |x|^{d-1}$.
\end {lemma}

This is one of many results that have been recently derived about connections in half-spaces, but we state here just this particular upper bound (where one of the two points lies on the boundary of the half-space). Indeed, the papers \cite {CH,CHS,panis_sharp_2025} all contain more precise up-to-constant estimates for connection probabilities in half-spaces. The latter paper \cite {panis_sharp_2025} in fact provides a uniform up-to-constant expression for $P [ x \lra_H y ]$ valid for all $x,y$.

\begin{figure}[h]
  \centering
  \includegraphics[width=\textwidth]{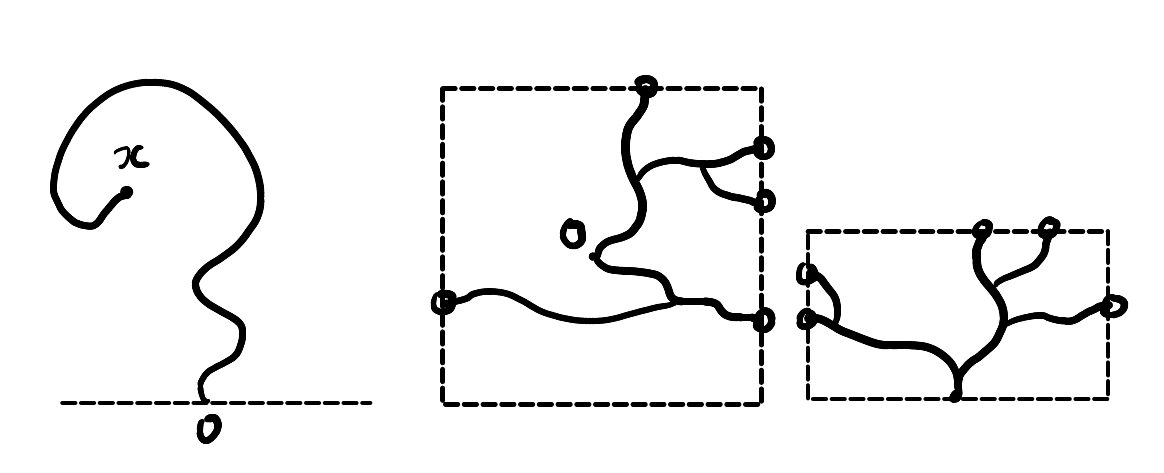}
  \caption{\label{w3}Events and quantities involved in Lemmas \ref {L4}, \ref {L5} and \ref {L6} respectively}
  \end{figure}

\subsection {Expected number of touching point for boxes}
Let us now state a closely related result in the whole space that will be useful in our setting, and that is also in a similar spirit to Principle ${\mathcal P}_1$.
Let $y$ be the point $(2n, 0, \ldots, 0)$ in $\Z^d$. When  $0 \lra y$, there necessarily exists $t \in \partial \Lambda_n$ such that $(0 \lra_{\Lambda_n} t) \circ (t \lra y)$. This immediately implies the bound
$$ P [ 0 \lra y] \le \sum_{t \in \partial \Lambda_n} P[0 \lra_{\Lambda_n} t ] P [ t \lra y], $$
and since in this inequality, both $P[ 0 \lra y] \asymp P[t \lra y ] \asymp 1/ |y|^{d-2}$, we see that the quantity
$$\Sigma_n := \sum_{t \in \partial \Lambda_n} P [ 0 \lra_{\Lambda_n} t] = E [\# \{ t \in \partial \Lambda_n \ : \  0 \lra_{\Lambda_n} t \} ]$$
is bounded from below by a constant that is independent of $n$.
The very useful feature here is that it is also bounded from above:
\begin {lemma}[\cite {CH}]
\label {L5}
There exists a constant $C$ such that $\Sigma_n \le C$ for all $n$ (hence $\Sigma_n \asymp 1$).
\end {lemma}
There are numerous ways to ``explain'' this result -- the perspective may differ a little depending on which way one approached the derivation of the two-point function estimate. Indeed, this result is in fact quite closely related to the two-point estimate itself. So, it is a little bit a matter of chicken and egg which one is derived first. In any case, this result has appeared explicitly in a number of recent papers, in particular the papers that derive restricted two-point functions in half-spaces \cite {CH,CHS,panis_sharp_2025}. Let us here however very briefly explain that it can be viewed as a direct consequence of these half-space estimates to stress the close relation between Lemmas \ref {L5} and \ref {L4}:

By reverting the roles of $t$ and $y$ and replacing them formally by $0$ and $x$ in this expression (and noting that the box is contained in a half-space, we obtain a uniform upper bound
$$ P [ 0 \lra_{\Lambda_n} t ] \le \frac {1} { n^{d-1}} $$
that is uniform over $n$ and $t \in \partial \Lambda_n$. Summing over the $O(n^{d-1})$ points in $\partial \Lambda_n$ then indeed leads to Lemma \ref {L5}.

It is worth noticing that this use of theses flows are not unrelated to the recent very effective approach to sharpness of the phase transition of percolation by
Duminil-Copin and Tassion (via their function $\phi(S)$) in \cite {DCT}.

\subsection {Expected number of touching points for half-boxes}
Let us state an analogous result to Lemma \ref {L5} in half-spaces:
Let us define
$$ \Sigma_n' := \sum_{t \in H \cap \partial \Lambda_n} P [ 0 \lra_{H \cap \Lambda_n} t]= E [\# \{ t \in \partial \Lambda_n \cap H \ : \  0 \lra_{\Lambda_n \cap H} t \} ].$$

\begin {lemma}[\cite{CH}]
\label {L6}
There exists a constant $C$ such that $\Sigma_n' \le C/n$ for all $n$.
\end {lemma}
The result appears as Lemma 26 in \cite {CH} and it is also part of Proposition 2.2 in \cite {panis_sharp_2025}.
This quantity is somehow related to $P [ 0 \lra_H x]$ in a similar way as $\Sigma_n$ was related to $P [ 0 \lra x]$. One can in particular deduce the upper bound for $P[0 \lra_H x]$ from this lemma (and this is indeed the order in which \cite {CH} proceeded to obtain Lemma \ref {L1}), noting that if $|x| > 2n$,
$$ P [ 0 \lra_H x ] \le  \sum_{t \in H \cap \partial \Lambda_n} P [ (0 \lra_{H \cap \Lambda_n} t) \circ (t \lra x)] \asymp \frac{1}{n^{d-2}} \times \Sigma_n'.$$

\subsection {Connections within a box}
A further remark (also pointed out in \cite {CH}) goes
as follows:
\begin {lemma}[\cite {CH}]
\label{L:restrict}
For any large enough fixed  $M>1$, one has
$$ P [ x \lra_{\Lambda_{MN}} y] \asymp P [ x \lra y ] \asymp \frac 1 { |x-y|^{d-2}}$$
for all $x \not= y$ in $\Lambda_{N}$.
\end {lemma}
This result holds in fact for any $M >1$ (see \cite {CH}), but we won't need this here and the direct proof below works easily when $M$ is chosen large enough, so we state only this case.
\begin{figure}[h]
  \centering
  \includegraphics[width=
  \textwidth]{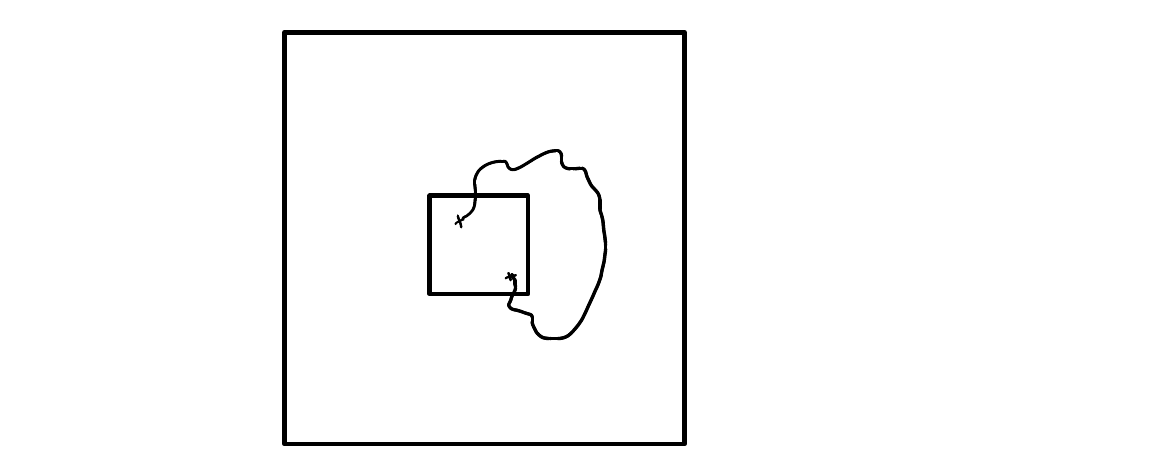}
  \caption{\label{fig:rest}Connections within a box $\Lambda_{MN}$ between points in $\Lambda_N$ (Lemma \ref {L:restrict})}
  \end{figure}
\begin {proof}
The upper bound is obvious since $P [ x \lra_{\Lambda_{MN}} y] \le P [ x \lra y] $.
For the lower bound, just note (using the half-space lemma in the case where one point is on the boundary of the half-space) that
\begin {eqnarray*}
\lefteqn {P [ (x \lra y) \setminus (x \lra_{\Lambda_{MN}} y) ]
 \le \sum_{t \in \partial \Lambda_{MN} } P [ (x \lra_{\Lambda_{MN}} t) \circ (t \lra y)]}
 \\ &&  \le C \frac {(MN)^{d-1}}{((M-1)N)^{d-1} ((M-1)N)^{d-2}} \le \frac {C'}{M^{d-2}} P [ x\lra y] \end {eqnarray*}
(where $C$ and $C'$ do not depend on $M$) which is smaller than $P[x \lra y] / 2$ provided that $M$ is large enough. Hence,
$ P [ x \lra_{\Lambda_{MN}} y] \ge P [ x \lra y] /2$
which allows us to conclude.
\end {proof}

\subsection {Some general comments about percolation in high dimensions}
\label {Scomment}

We now make a little pause, and briefly discuss some of the general principles that lie behind most of the proofs in the study of
high-dimension critical percolation (with the notable  exception of the lace-expansion technique used to derive the two point estimate (\ref {tpe})) -- this section corresponds to what was Section 2 in the first version of the present paper on arXiv. In particular, let us stress the principle that  the upper bounds that are coming from the use of the BK inequality to decompose probabilities end up being sharp. This is best seen by listing examples:

\begin {enumerate}
 \item
The first simple instance is that
 $$ P [ (x \lra t ) \circ (t \lra y) ] \asymp P [ x \lra t ] P [ t \lra y].$$
 We will describe this is Section \ref {Sdc} in the  case where $x=y$ -- then, $(x \lra t ) \circ (t \lra x)$  corresponds to the existence of a loop-cluster going through $x$ and $t$).
The upper bound is just the BK inequality, and one strategy to prove the lower bound will be explained in that Section.
 \item
The next example is given by multiple-point connectivity functions, where trifurcation diagrams show up (this is very much in the spirit of Aizenman's paper \cite {Ai}). Clearly, the BK inequality shows that (here $t$ denotes the first point at which an open path from $z$ to $x$ meets a self-avoiding open path from $x$ to $y$):
$$ P [ x \lra y \lra z ] \le \sum_{t} P [ (x \lra t) \circ (t \lra y) \circ (t \lra z)] \le \sum_P P [ x \lra t] P[ y \lra t] P[ z \lra t].$$
As it turns out, the right-hand side is in fact bounded by a constant times $P [ x \lra y \lra z ]$. This roughly corresponds to the idea that when $x, y, z$ are in the same cluster, then the number $T$ of points $t$ for which  $(x \lra t) \circ (t \lra y) \circ (t \lra z)$ holds is small (i.e., tight). To shows, one natural possible strategy is to show that $E[ T^2]$ is bounded by a constant times $E[ T]$. So, one needs a sharp lower bound on $P[T]$ (that can be obtained in a similar manner to that of the first example) and then to estimate the second moment of $T$, where one can  use a diagrammatic expansion that can usually be controlled easily.
 \item
 The more involved but crucially important case of flows out of boxes has already been mentioned; the upper bounds of the type
 $$ P [ x \lra y ] \le \sum_{t \in \partial \Lambda} P[(x \lra_\Lambda t)\circ (t \lra y)] \le \sum_{t \in \partial \Lambda} P [ x \lra_\Lambda t] P [ t \lra y]$$
 when $x \in \Lambda$ and $y \notin \Lambda$
 are sharp as well (i.e., the right hand side is in fact bounded by a constant times the left-hand side). This roughly corresponds to the idea that when $x \lra y$, then the number $T$  of points $t$ for which $(x \lra_\Lambda t)\circ (t \lra y)$ holds is typically very small (i.e. tight).
 Things are a bit trickier to explain the proof
 as we have the chicken and egg issue now -- the previous diagrammatic estimates require the uniform half-space upper bounds that are themselves related to the result itself. But to control these diagrammatic estimates, one just needs to show convergence of some sums, so some room is available in the upper bounds that one uses (this is partially  what the aforementioned papers are about).
\end {enumerate}
A recently posted paper by Panis and Schapira \cite {panisschapira2} pushes those ideas further -- in particular showing that the constants in the last bound can be chosen uniformly with respect to the set $\Lambda$ (as in \cite {DCP}, they call this a reversed Simon-Lieb inequality, in reference to the corresponding correlation inequalities for spin systems or $\phi^4$ models in \cite {Simon,Lieb}), leading to new derivations of various results, such as a new proof of the Kozma-Nachmias \cite {KN2} one-arm probability  estimate (in $\Z^d$)
\begin {equation}
 \label {onearm}
P [ 0 \lra \partial \Lambda_n ] \asymp \frac 1 {n^2}.
\end {equation}

In this paper, we will mostly use upper bounds, but in order to prove the existence part of  Proposition \ref {mainprop}, we will use the general ideas that we have just outlined.

\section {Self-avoiding connections between nearby points exiting a large box}

We now derive the lemmas that will then be used to derive the propositions stated in the introduction,

\subsection {Connections exiting a large box}

We now focus on a slightly different type of connections, where the open self-avoiding path between $x$ and $y$ has to go out of a big set. More specifically, we will look at:
\begin {itemize}
 \item
``Long'' connection events: When $x$ and $y$ are two points in a set $\Lambda$, we denote by  $x \lra^\Lambda y$ the event that there exists an open self-avoiding path joining $x$ and $y$ that exits $\Lambda$.
\item
``Long restricted'' connections: If $x$ and $y$ are two points in $\Lambda_1$ and $\Lambda_1 \subset \Lambda_2$, we denote by $x \lra^{\Lambda_1}_{\Lambda_2} y$ the event that there exists a self-avoiding open path joining $x$ and $y$ that stays in $\Lambda_2$ but exits $\Lambda_1$ (mind that this is not the same as the intersection of $x \lra^{\Lambda_1} y$ with $x \lra_{\Lambda_2} y$ since we now require both events to be realized by the same open path).
\end {itemize}

The following simple lemma will be useful:

\begin{lemma}\label{L:horseshoe}\label{L-horseshoe}
    There exists a constant $C$ such that for all $x\in \Lambda_{N}$ and $N\geq 1$,
  $ P [ 0 \lra^{\Lambda_N} x ] \le C/N^{d-2}$.
\end{lemma}
\begin {proof}
This is a direct consequence of Lemma \ref {L5}. Indeed, it clearly suffices to consider the case where $x \in \Lambda_{N/2}$ (otherwise the upper bound $P[ 0 \lra x] $ does the job), and then
$$   P [ 0 \lra^{\Lambda_N} x ] \le \sum_{t \in \partial \Lambda_N} P [ 0 \lra_{\Lambda_N} t] P [ t \lra x ] \le C N^{d-1} \asymp \frac {1}{N^{d-2}} \sum_{t \in \partial \Lambda_N} P [ 0 \lra_{\Lambda_N} t] \le \frac {C'}{N^{d-2}} $$
for some constant $C'$ by Lemma \ref {L5}.
\end {proof}
Note that when $x=0$, this indicates that the probability that the origin (or any other fixed point) is on a loop-cluster of diameter at least $N$ is comparable (up to multiplicative constants) to
$ 1 / N^{d-2}$. In particular, this immediately shows that the expected number of points in $\Lambda_n$ that are in a loop-cluster of diameter at least $an$ is bounded by a constant times $n^2$.

The astute reader might actually guess (from the general principles outlined in Section \ref {S2}  that the bound in the lemma is in fact sharp, i.e., that
$$ P [ 0 \lra^{\Lambda_N} x ] \asymp \frac 1 {N^{d-2}}$$
for $x \in \Lambda_{N/2}$. This is indeed the case and can in fact be proved using the same ideas that we outlined, but for the purpose of the present paper, this upper bound will be sufficient.

\begin{figure}[h]
  \centering
  \includegraphics[width=\textwidth]{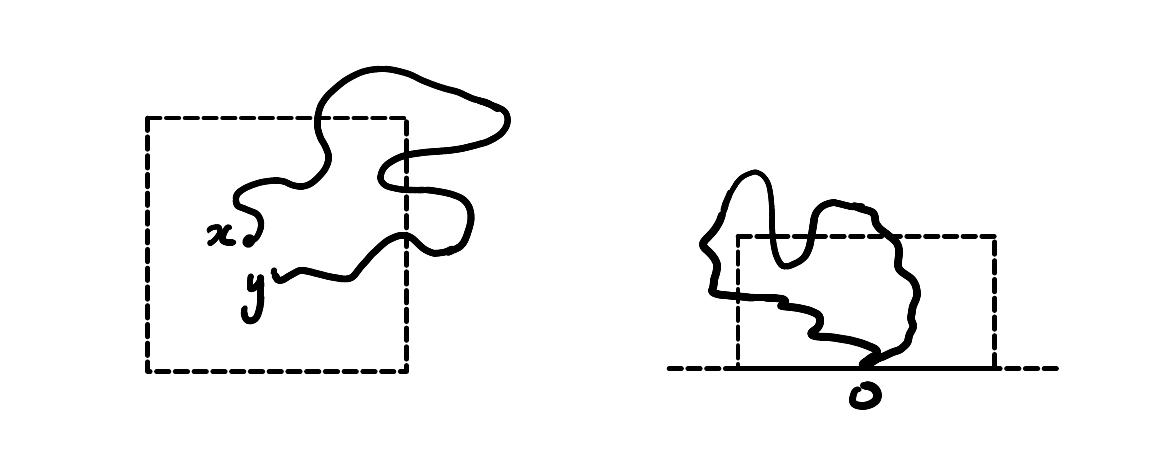}
  \caption{\label{w4}Events $x \lra^{\Lambda_N} y$ and $0 \lra^{\Lambda_N}_H 0$ involved in Lemmas \ref {L:horseshoe} and \ref {L:restrict2}}
  \end{figure}

\subsection {A half-space version and its consequence}

\label {Supperbound}

In this section, we derive an analogous result for long connections exiting a box in a half-space, and then deduce the tightness of the number of large loop-clusters.

A large loop-cluster will have a point with lowest first coordinate. So, seen from such a point $x=(x_1, \ldots, x_d)$, the large loop will be contained in the half-space consisting of all points with first coordinate greater or equal to $x$. So, it is natural to want to estimate the corresponding probabilities for each given fixed $x$. In other words, if $x$ is the origin and $H$ denote the half-space consisting of all points with non-negative first coordinate:
\begin {lemma}
\label{L:restrict2}
For some constant $C$,  $ P [ 0 \lra_{H}^{\Lambda_n} 0] \le  C/n^d$.
\end {lemma}

The astute reader might again already infer that this upper bound is sharp, i.e., that
$$P [ 0 \lra_{H}^{\Lambda_n} 0] \asymp \frac 1 {n^d}$$
which is indeed true (and could be proved using similar techniques as in this paper).

\begin {proof}
The proof  goes along the now familiar lines as in the whole space, except that we now use the bound on $E [ \Sigma_n']$ instead of the bound on $E[ \Sigma_n]$: We first use the BK inequality and Lemma \ref {L6} to see that
$$ P [ 0 \lra_{H}^{\Lambda_{n}} 0] \le \sum_{t \in \partial \Lambda_{n}} P [ 0 \lra_{\Lambda_n \cap H} t] P[t \lra_H 0 ] \le \frac {C}{n^{d-1}}
\sum_{t \in \partial \Lambda_n} P [ 0 \lra_{\Lambda_n \cap H} t].$$
But this last sum is bounded by a constant times $1/n$ by Lemma \ref {L6}, which concludes the proof.
\end {proof}

\begin {proof}[Proof of tightness in Proposition \ref {mainprop}]
This lemma immediately provides an upper bound for the expected number of  loop-clusters of $L^\infty$-diameter greater than $a N$ in $\Lambda_N$.
Indeed, any self-avoiding loop of $L^\infty$-diameter at least $a N$, when seen from one of its points $x$ with lowest first coordinate will be a self-avoiding loop of diameter at least $a N$ in the half-space ``above'' $x$ (i.e., the set of points with first coordinates at least equal to that of $x$). For each $x$, the probability of this event is bounded by $C(a)/ N^{d}$ by Lemma \ref {L:restrict2}, so that
 summing over all $O(N^d)$ points $x$ in $\Lambda_N$ shows that the expected number of self-avoiding loops of diameter greater than $a N$ is bounded by some finite constant.
In other words, the tightness part of Proposition \ref {mainprop} follows.
Note that this does however not yet prove the actual existence of macroscopic loops.
\end {proof}

\section {Proof of Proposition \ref {main2} using Lemma \ref {L-horseshoe}}

In this section, we will use the upper bounds of Lemma \ref {L-horseshoe} to deduce a couple of features of loop-clusters. The idea will essentially to upper-bound some probabilities or expectations in order to deduce that the probability of certain events do vanish in the large $N$ limit. This will in particular imply Proposition \ref {main2}.

\subsection {No two disjoint big loops in the same cluster}
Consider percolation in $\Lambda_N$. Let $a<1$ and $k$, and let $E(N,k, k')$ denote the event that there exists a percolation cluster that contains two disjoint self-avoiding loops, one of $L^\infty$ diameter greater than $k$ and one of $L^\infty$ diameter greater than $k'$.

Note that if $E(N,k, k')$ holds, then one can then find an open path joining the two disjoint loops, so that there exist two points $z$ and $z'$ such that
$$ (z \lra^{z +\Lambda_{k/2}} z) \circ (z \lra z') \circ (z' \lra^{z'+ \Lambda_{k'/2}} z').$$
By the BK inequality, Lemma \ref {L:horseshoe} and summing over $z$ and $z'$, we immediately get an upper bound of a constant times
$$
 \frac {1}{k^{d-2} (k')^{d-2}} \times N^d \times N^2.
$$
This immediately implies that some configurations will not occur in the $N \to \infty$ limit:
In particular, if we fix some small $a$ and choose $k' = aN$, we see that
$ P [ E(N, aN, k )]$ is bounded by a constant times $N^4 / k^{d-2}$.
So, this for instance implies that for any $m_N \to \infty$, if we set $u(N)= m_n N^{4/ (d-2)}$  (mind of course that $4/(d-2) \le 4/5< 1$ since $d \ge 7$):
\begin {lemma}\label{L:no-glasses}
For fixed $a$, the probability that some loop-cluster in $\Lambda_N$ contains two disjoint self-avoiding loops, one with diameter at least $aN$ and one with diameter at least $u(N)$, goes to $0$ as $N \to \infty$.
\end {lemma}

\subsection {No two essentially different big loops intersect}
Consider percolation in $\Lambda_N$. Let   $F(N,k, k')$ denote the event that there exist two points $x$ and $y$ such that the following three events occur disjointly:
$( x \leftrightarrow y )$, $( x \leftrightarrow^{x + \Lambda_k} y )$ and $( x \leftrightarrow^{x + \Lambda_{k'}} y)$.

In particular, if there is a point $z$ in the loop-cluster that is at distance greater than $k'$ from a loop of diameter $2k$, this event will necessarily hold (indeed, the first two connections then correspond to the loop of diameter at least $2k$, and the third connection correspond to the excursion containing $z$ of the second loop away from the first one). It more generally says that any open self-avoiding excursion away from the first loop has diameter at most $k'$,

Again, the probability of $F(N,k,k')$ is immediately bounded by combining the BK inequality, Lemma \ref {L:horseshoe} and summing over $x$ and $y$
[when $x$ and $y$ are at distance smaller than $k$ from each other, we use Lemma \ref {L:horseshoe} to bound $P [  x \leftrightarrow^{x + \Lambda_k} y]$ and if $x$ and $y$ are at distance greater than $k$, we can use the usual two-point function bound $\asymp 1/ |x-y|^{d-2}$].

For instance, if we choose $k'=aN$ for some fixed $a$, then we get that
$P [ F(N, aN, k)]$ is bounded by a constant times $N^{2}k^{4-d}$. In particular, $$ P [  F(N, aN, m_n N^{2/ d-4}) ] \to 0$$ as $N \to \infty$ if $m_N$ is chosen in such a way that $u_N \to \infty$.

This  implies that for any $m_N \to \infty$ if we now set $l(N):= m_N N^{2/(d-4)}$ (mind of course that $2/(d-4)<4/(d-2)$ is also smaller than $1$ since $d \ge 7$):
\begin {lemma}
\label {L:nothetashape}
For fixed $a$, the probability that some loop-cluster in $\Lambda_N$ contains two self-avoiding loops of diameter at least $aN$ that are at Hausdorff distance greater than $l(N)$ from each other or from $\gamma_1 \cap \gamma_2$ tends to $0$ as $N \to \infty$.
\end {lemma}

\begin{figure}[h]
  \centering
  \includegraphics[width=
  \textwidth]{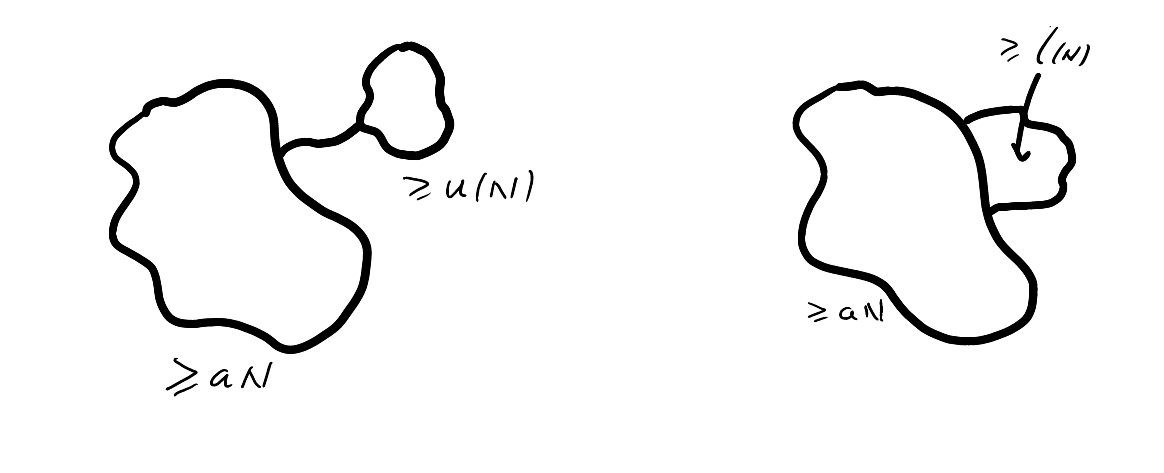}
  \caption{\label{fig:4}Schematic depiction of the events involved in Lemmas \ref {L:no-glasses} and \ref {L:nothetashape}}
  \end{figure}

 \begin {proof} By Lemma \ref{L:no-glasses}, any two large self-avoiding loops $\gamma_1$ and $\gamma_2$ (of $L^\infty$ diameter greater than $aN$) that are in the same percolation cluster do necessarily intersect (with probability that goes to $1$ as $N \to \infty$). They are therefore part of the same loop-cluster with probability $1 -o (1)$. Suppose that $z$ is a point on $\gamma_1$ that is at distance at least $l(N)$ from $\gamma_1 \cap \gamma_2$. Let $z$ be such a point on $\gamma_1$, we can then follow $\gamma_1$ in both directions starting from $z$ up to the first points $x$ and $y$ at which it intersects $\gamma_2$ -- these points will both be at distance at least $l(N)$ from $z$.  These two points are therefore joined by three disjoint self-avoiding open paths: The two portions of $\gamma_2$ that join them (one of which having diameter at least $aN/2$), and the portion of $\gamma_1$ that we have just defined. But the probability of the existence of such a configuration goes to $0$ by our previous bound on the probabilities of the events $F$.
 So the probability that the Hausdorff distance between $\gamma_1$ and $\gamma_1 \cap \gamma_2$ is greater than $l(N)$ goes to $0$, and since the same is true when interchanging the roles of $\gamma_1$ and $\gamma_2$, the lemma follows.
 \end {proof}

\subsection {Self-avoiding on sites versus  self-avoiding on edges}

Lemma \ref {L:nothetashape} ensures that when $m_N \to \infty$,  with a probability that tends to $1$ as $N \to \infty$, any loop-cluster of diameter greater than $aN$ will in fact contain a loop of diameter at least $u_N := aN - m_n  N^{2/ d-4}$ that is self-avoiding on sites (and not just on edges). Indeed, if one chooses two points $x$ and $y$ at distance greater than $aN$ on the edge-self-avoiding open loop, then one can consider the loop-erasure of both edge-self-avoiding paths joining $x$ and $y$ in the loop, and one gets two site-self-avoiding paths $\gamma_1$ and $\gamma_2$ joining $x$ and $y$ using different edges. If these two paths meet at a point $z$ that is at distance greater than $u_N$ from $x$ and from $y$, then the union of the two paths will contain two edge-self-avoiding loops of diameter greater than $aN/2$ and $u_N /2$ respectively. By Lemma \ref {L:nothetashape}, this happens with probability that that goes to $0$ as $N \to \infty$.
Hence, with probability that goes to $1$, no such point $z$ exists (and the union of such $\gamma_1$ and $\gamma_2$ will contain a site-self-avoiding loop of diameter at least $aN - u_N$).

Alongside Lemmas \ref {L:horseshoe} and \ref {L:nothetashape}, this observation completes the proof of Proposition \ref {main2}.

\section {Probabilities of double-connection events}
\label {Sdc}
Suppose that $x$ and $y$ are two distant points, and let us estimate the probability that $x$ and $y$ are in the same loop-cluster.
The BK inequality and (\ref{tpe}) immediately give the upper bound
$$ P [ ( x \leftrightarrow y ) \circ ( y \leftrightarrow x ) ] \le P [ x \leftrightarrow y]^2 \asymp \frac 1 {|x-y|^{2d-4}}.$$
Let us now explain how to derive the matching lower bound, i.e., to show that
\begin {lemma}
\label {lemma1}
\label {L1}
One has (for all $x \not= y$)
$$ P [ ( x \leftrightarrow y ) \circ ( y \leftrightarrow x ) ]  \asymp \frac {1} {|x-y|^{2d-4}}.$$
\end {lemma}

\begin {proof}
Let us consider two independent critical percolation configurations $\omega$ and $\omega'$. By (\ref{tpe}), the probability of the event $E$ that $x \leftrightarrow y$ holds for both configurations satisfies $P [E] \asymp 1/|x-y|^{2d-4}$. Furthermore, the expected number $M$ of points $z$ such that (a) the event $E$ holds, (b) $z$ is in the cluster of $x$ for the first percolation $\omega$, and (c) $( x \leftrightarrow z ) \circ ( z \leftrightarrow y )$ for the second configuration $\omega'$ is  bounded by a constant times $P[E]$. Let us now detail how to show this:
For a given $z$, for the event involving $\omega$ to hold, it means that for some $t$, $ (x \lra t) \circ (t \lra y) \circ (t \lra z)$.
Using the BK inequality, and then summing over $t$ and then $z$ gives an upper bound for $E[M]$ of the type (we use the convention that $1/0 = 1$ is such sums)
$$ \asymp \sum_{z, t}  \frac {1} {|x-t|^{d-2}  |z-t|^{d-2} |y-t|^{d-2} |x-z|^{d-2}  |z-y|^{d-2} } .$$
To bound this sum, we will treat separately the cases where $z$ and/or $t$ are close to $\{ x, y\}$ or not.
Recall the Aizenman-Newman \cite {AN} triangle condition is satisfied, i.e., that
$$\sum_{u,v} \frac 1 { |u|^{d-2} |v|^{d-2} |u-v|^{d-2}} < \infty.$$
This immediately implies that  the partial sum over all $t$ and $z$ that are both at distance at most $3|x-y|/ 4$ of $x$ will be bounded by a constant times $(x-y)^{4-2d}$ (just note that then $|t-y| \ge |x-y|/2$  and $|z-y| \ge |x-y| /2$ and then apply the triangle condition). The same clearly applies when both $t$ and $z$ are both at distance at least $3|x-y|/ 4$ from $y$.
If  we sum over $x$ and $t$ such that $|x-t|< 2|x-y|/3$ and $|z-y|< 2|x-y|/3$, we get an upper bound
$$ \asymp \frac {1}{|x-y|^{3d-6}} \Bigl(\sum_{z , \ |z-y| < 2|x-y|/3}  1/|z-y|^{d-2}\Bigr)^2 \asymp  \frac { |x-y|^4}{|x-y|^{3d-6}}  =  \frac {1}{|x-y|^{2d-4 + d-6}}.$$
And finally, if one of the points $x$ or $t$ is at distance greater than $|x-y|/4$ of both $x$ and $y$ is treated easily. So, we conclude that the conditional expectation of $M$ given $E$ is indeed finite.

The same argument shows actually that when $K$ is fixed and large enough, if one restricts the sum over $z$ to the points that are at distance at least $K$ from both $x$ and $y$ we get that the expected number (conditionally on $E$) of such points is bounded by
$P[E]/2$ (for any $N$). In particular, it implies that the probability that $E$ holds and that the set of points at distance at least $K$ from $x$ and $y$ that satisfy (a), (b) and (c) is empty is at least $P[E]/2$. In this case, $\omega$ has an open path from $x$ to $y$ and $\omega'$ has an open path joining the boxes of size $K$ around $x$ to the box of size $K$ around $y$ and this path does not intersect the cluster containing $x$ for $\omega$.

By then resampling both configurations $\omega$ and $\omega'$ in the $K$-neighborhoods of $x$ and $y$ (noting that there are only finitely many possible resampling options and that at least one will do the job), one then concludes that if we denote by $E'$  the event that $E$ holds and that there exists an open path of edges for $\omega'$ that joins a given neighbour $x'$ of $x$ to a given neighbour $y'$ of $y$ that stays at distance at least one from the open cluster of $x$ for $\omega$, then $P [ E'] \ge cP[E]$ for some positive constant $c$.

By now first revealing the cluster containing $x$ for $\omega$, and then revealing the percolation status of the remaining edges (for which one can equivalently use the values taken by $\omega'$ instead of those of $\omega$ since the outcome will have the same law), we therefore see that (if $x'$ is the chosen neighboring point of $x$ and $y'$ is the chosen neighboring point of $y$), then
$$ P [ ( x \leftrightarrow y ) \circ ( x' \leftrightarrow y' ) ] \ge c P[ x \leftrightarrow y ]^2 .$$
Finally, by resampling the status of the edge between $x$ and $x'$ and the edge between $y$ and $y'$, we conclude that
$$ P [ ( x \leftrightarrow y ) \circ ( x \leftrightarrow y ) ] \ge c' P[ x \leftrightarrow y ]^2 $$
for some constant $c'$, which provides the matching lower bound.
\end {proof}

\begin {remark}
The very same proof can be adapted to derive the following statements involving linear chains of connections:
$$ P [ ( x \leftrightarrow t ) \circ ( t \leftrightarrow y ) ] \asymp P[ x \leftrightarrow t ] P [ t \leftrightarrow y]$$
and
$$ P [ ( x_1 \leftrightarrow x_2 ) \circ ( x_2 \leftrightarrow x_3 ) \circ \cdots \circ ( x_{n-1} \leftrightarrow x_n ) ]
\asymp \prod_{j=1}^{n-1} P [ x_j \leftrightarrow x_{j+1} ].$$
Similarly, for circular connections (i.e., loops),
$$ P [ ( x_1 \leftrightarrow x_2 ) \circ ( x_2 \leftrightarrow x_3 ) \circ \cdots \circ ( x_{n-1} \leftrightarrow x_n ) \circ (x_n \lra x_1)  ]
\asymp P[ x_n \lra x_1] \times \prod_{j=1}^{n-1} P [ x_j \leftrightarrow x_{j+1} ] .$$
The proof also works for any finite connected diagrams of connections, for instance
to show that
$$P [ (x_1 \lra y) \circ (x_2 \lra y) \circ (x_3 \lra  y)] \asymp \prod_{j=1}^3 P[ x_j \lra y].$$
\end {remark}

\begin {remark}
The diagrammatic sums in the proof are one instance of the Aizenman-Newman  idea that in such sums, one can basically perform a triangle star transformation for all inner triangles in the diagram (and only compensate by a multiplicative constant) -- which is another of the general principles that have been used in numerous occasions in this field.
\end {remark}

The same argument as in the proof of Lemma \ref {L:restrict} then gives the similar result for double connections restricted to a box:
\begin {lemma}
\label {lemma2}
\label {L2}
For any  large enough $M>1$, one has
$$ P [ (x \lra_{\Lambda_{MN}} y)\circ (x \lra_{\Lambda_{MN}} y)] \asymp P [ x \lra y ]^2 \asymp \frac 1 { |x-y|^{2d-4}}$$
for all $x \not= y$ in $\Lambda_{N}$.
\end {lemma}
\begin {proof}
One has
\begin {eqnarray*}
  \lefteqn{P [ (x \lra y)\circ (x \lra y) \setminus (x \lra_{\Lambda_{MN}} y) \circ (x \lra_{\Lambda_{MN}} y) ]}\\
 &\le&  2 \sum_{t \in \partial \Lambda_{MN} } P [ (x \lra_{\Lambda_{MN}} t) \circ (t \lra y) \circ (y \lra x)]\\ &  \le&  \frac {C (MN)^{d-1}}{((M-1)N)^{d-1} ((M-1)N)^{d-2} |x-y|^{d-2}} \le \frac {C'}{M^{d-2}}\PP(x\lra y)^2\end {eqnarray*}
(where $C$ and $C'$ do not depend on $M$) which is smaller than $P[(x \lra y) \circ (x \lra y)] / 2$ provided that $M$ has been chosen to be large enough.
\end {proof}

\section {Existence of large loop-clusters via a second moment bound}

So far, we have not yet proved that large loop-clusters do exist in the large $N$ limit.
There are several possible ways to proceed. One natural option is to follow the idea mentioned in the introduction, namely to use the
framework of Section \ref {Supperbound}, and to define the number $S$ of points with lowest first coordinate on a self-avoiding open loop of diameter at least $aN$. The strategy is then to show that $E[S]$ is bounded from below by a constant, and that $E[S^2]$ is bounded from above by a constant (then, with Cauchy-Schwarz, this implies that $P[ S \not= 0]$ is bounded from below). We will outline in Section \ref {Sbrief} how to implement this strategy building on the uniform version of the half-space two-point functions from \cite {panis_sharp_2025}, but we opt here for a detailed presentation of a slightly different approach, that turns out to be a little simpler: Instead of focusing on $S$, we will estimate the first and the second moment of the number $\pi$ of pairs of points $x$ and $y$ in $\Lambda_N$ that are at distance greater than $aN$ from each other and also in the same loop-cluster. This allows to avoid half-space estimates altogether.

\subsection {Via pairs of distant points}
So we now focus on the first and second moments of $\pi$.
By summing over all points $x$ and $y$ in $\Lambda_N$, it follows from our previous estimates that (provided $a$ has been chosen small enough)
$$E [ \pi ] \asymp N^4 .$$
Indeed, the upper bound follows immediately from Lemma \ref {L1}, and for the lower bound, we can use Lemma \ref {L2} (this is where we need $a$ to be small enough, i.e. somewhat smaller than $1/M$).

The goal is then to bound the second moment of $\pi$:
\begin {lemma}
\label {lemma3}
One has $E [ \pi^2 ] \le C N^8$, for some constant $C=C(a)$.
\end {lemma}

Let us first explain how to conclude the proof of Proposition \ref {mainprop} using Lemma \ref {lemma3}:
\begin {proof}[Proof of existence in Proposition \ref {mainprop}]
One uses the usual  Cauchy-Schwarz based argument. We have that
$$ P [ \pi \not = 0 ] \ge E[ \pi]^2 / E[ \pi^2] $$
and it follows from Lemmas \ref {lemma3} and \ref {lemma2} that the right hand side is bounded away from $0$.
But the event that $\pi \not= 0$ is exactly the event that there exists at least one loop-cluster of diameter at least $aN$.
So, to summarize (switching notations slightly, replacing $a$ and $N$ by $\eps$ and $m$), the probability that when restricting the percolation to $\Lambda_m$, provided $\eps$ has been chosen small enough (but fixed), the probability that there exists a loop-cluster of diameter at least $\eps m$ in $\Lambda_m$ is bounded from below by a positive constant $c(\eps)$ independently of $m$.

To conclude the proof of Proposition \ref {mainprop}, we need to show that the numbers $n(a,N)$ of loop-clusters of diameter at least $aN$ will actually be very big when $a$ is chosen sufficiently  small. To see this, we can proceed as follows: Let $\eta>0$ and $k$ be fixed.
We first choose any small $\eps$ as above, with $c=c( \eps) > 0$.
 When $l$ is chosen large enough, the probability that $l^d$ independent Bernoulli variables with parameter $c$ contain more than $k$ positive outcomes is greater than $1- \eta /2$. Let us fix such an $l$ and define $a:= \eps / (2l)$

For any large $N$ large, we can find $m > N / (2l)$ and  $l^d$ disjoint boxes of the same size as $\Lambda_m$ in $\Lambda_N$.  We then apply the previous result to each of them (noting that the percolation outcomes inside each of these boxes are independent), and see that with a probability at least $1- \eta/2$, there will be at least $k$ of these $l^d$ boxes that contains loop-clusters (for the percolation restricted to that box) of diameter at least $\eps m$ (which is greater than $\eps N / (2l) = aN$).

But we also know from Lemma \ref{L:no-glasses} that the probability that there exist in $\Lambda_{N}$ two disjoint loops of diameter at least $aN$ that are in the same loop-cluster goes to $0$ as $N \to \infty$. it is smaller than $\eta/2$ for large enough $N$. Hence, when $N$ is large enough, the probability that there exist at least $k$ disjoint loop-clusters of diameter greater than $a N $ in $\Lambda_{N}$ is greater than $1- \eta$, which concludes the proof of Proposition \ref {mainprop}.
\end {proof}

\begin {proof}[Proof of Lemma \ref {lemma3}]
We want to sum over $x$, $y$, $x'$ and $y'$ in $\Lambda_N$ with the constraints that $|x-y|>aN$ and $|x'-y'|>aN$ the probability of the event
$$  ((x \lra y) \circ (x \lra y)) \cap ((x' \lra y') \circ (x' \lra y')).$$
We can consider separately the two subcases of this event where the loops containing $\{ x, y\}$ and $\{ x', y'\}$ can be chosen to be disjoint, and where these loops have to intersect.  The probability of the former case is bounded by the BK inequality, and summing these contributions just leads to the upper bound $\asymp E[ \pi ]^2 \asymp N^8$ as desired.  It remains to focus on the latter case, which is less straightforward.

In this case, a loop $l'$ that contains $\{x', y'\}$ will intersect a  loop $l$ that contains $\{x, y\}$. By only keeping the excursions  of $l'$ away from $l$ that do actually contain $x'$ and $y'$ (there could be one or two excursions) i.e., by discarding the remaining parts of $l'$, we end up with ten type of configurations denoted by (0)-(9) in Figures \ref {g1}, \ref{g2}, \ref {g3} and \ref {g4} (in those sketches, the points $x$ and $y$ are the white dots, and the points $x'$ and $y'$ are the black dots) depending on the relative positions of the endpoints of these excursions. We can then further note that some configurations of connections imply somewhat simpler ones (one can for instance modify $l$ in such a way that it goes through one or both of $x'$ and $y'$ in cases (0), (5), (6) and (8)). So, one then ends up with the set of configurations (2')-(9') in the figures. The configuration (9') can be rephrased in to (9'') by exchanging the roles/relative positions of $x',y'$ with those of $x,y$, which in turn is the same as (0) with the colors reversed.

\begin{figure}[h]
  \centering
  \includegraphics[width=.9\textwidth]{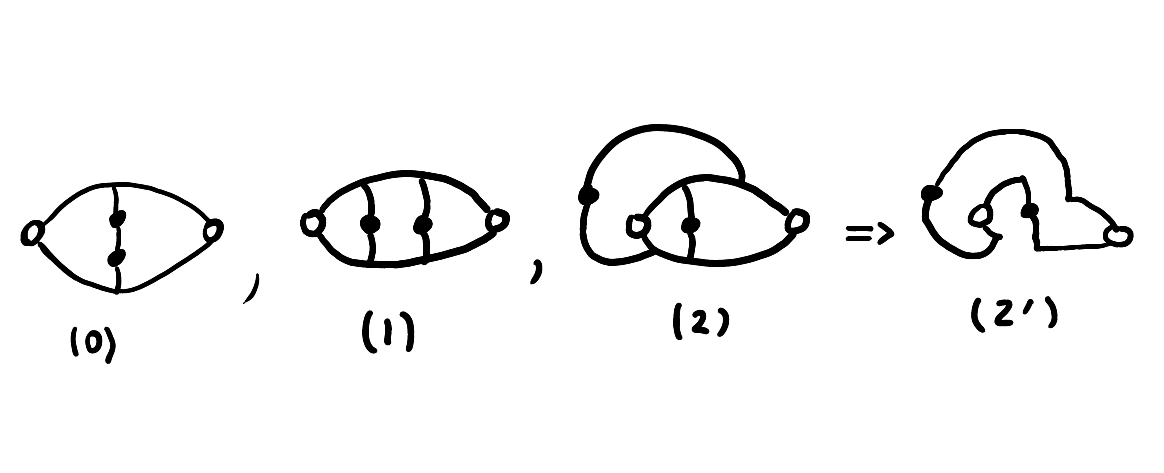}
  \caption{\label{g1}The cases when endpoints of each excursion of $l'$ land on opposite parts of $l$}
  \end{figure}

\begin{figure}[h]
  \centering
  \includegraphics[width=.6\textwidth]{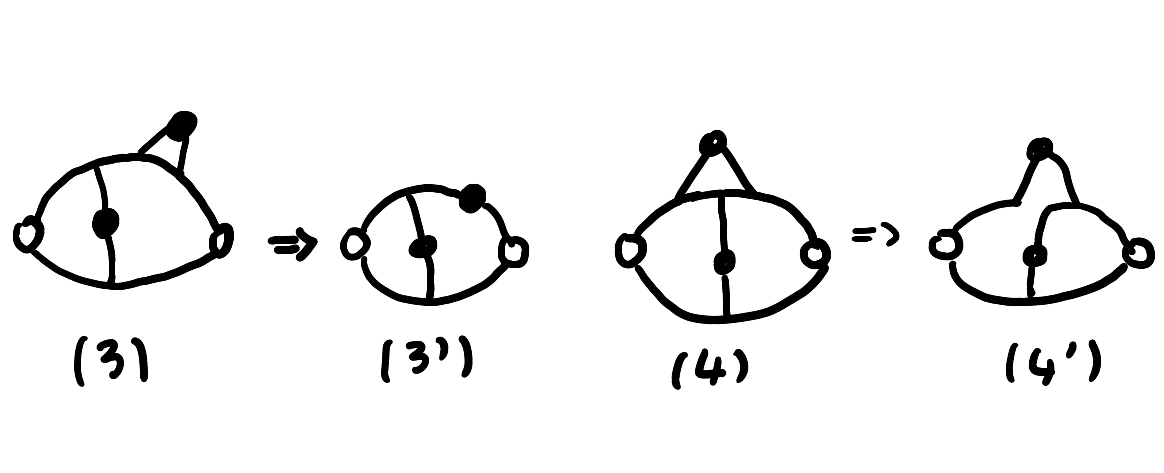}
  \caption{\label{g2}The cases when the endpoints of only one out of two excursions of $l'$ land on opposite parts of $l$}
  \end{figure}

\begin{figure}[h]
  \centering
  \includegraphics[width=.9\textwidth]{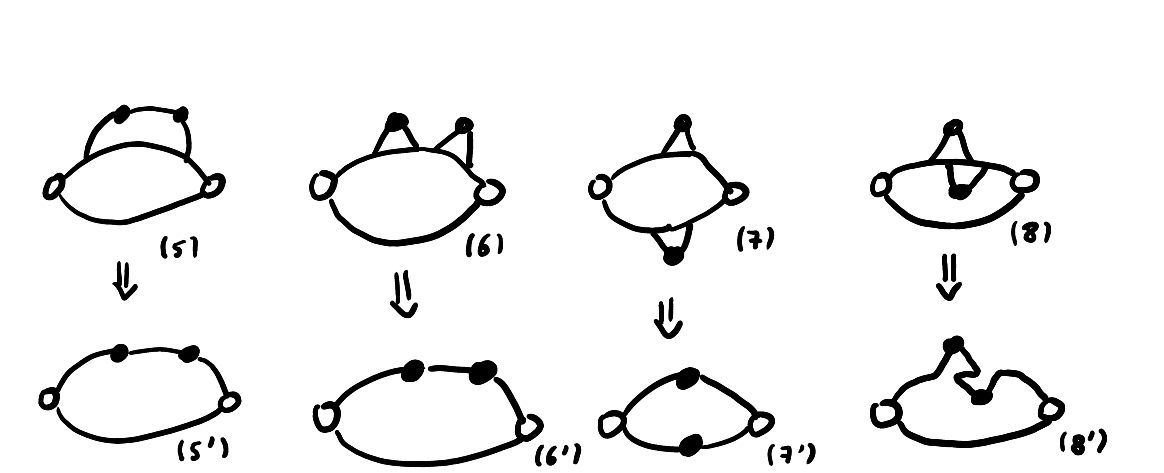}
  \caption{\label{g3}The first four remaining cases}
  \end{figure}

\begin{figure}[h]
  \centering
  \includegraphics[width=.6\textwidth]{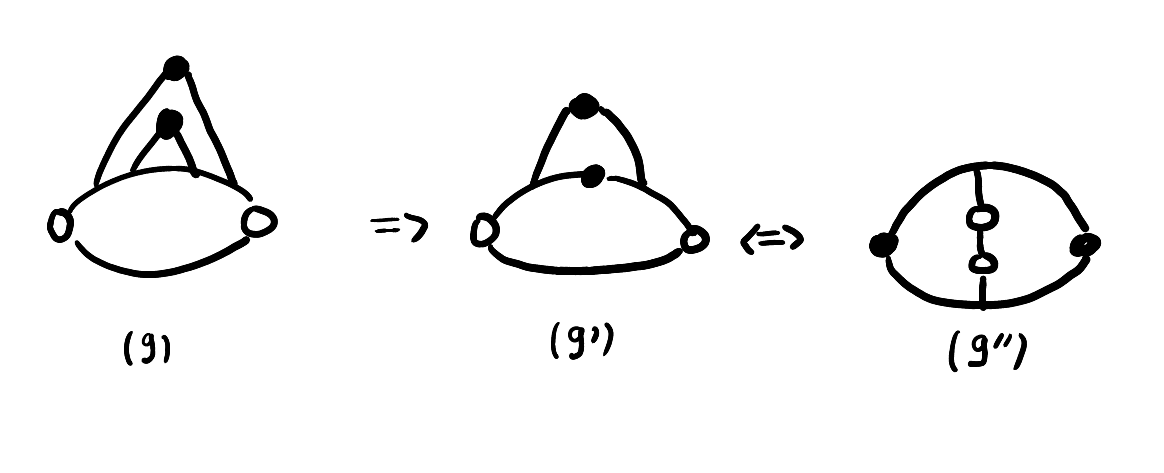}
  \caption{\label{g4}The final case}
  \end{figure}
Note that the probability of each individual configuration  has a simple upper bound obtained by the BK inequality and the upper bound on the two-point function.
The job is now to check that when one fixes $x$ and $y$ (i.e., that are at distance at least $aN$ from each other) and then sums these probabilities over all admissible $x'$ and $y'$ (also at distance at least $aN$ from each other), one has an upper bound of a constant times $N^{8-2d}$ (so that summing over all $O(N^{2d})$ options for $x$ and $y$ indeed provides the $O(N^8)$ upper bound).

This final useful input is given by the simple estimates (for each fixed $u \not= v$)
\begin {equation}
\label {final1}
\sum_{t} P [ (u \circ t ) \circ (t \circ v) ] \le C \frac {|u-v|^2}{|u-v|^{d-2}} \le C' \frac {N^2}{|u-v|^{d-2}}
\end {equation}
and
\begin {equation}
 \label {final2}
\sum_{t_1, t_2} P [ (u \lra t_1 ) \circ (t_1 \lra t_2) \circ (t_2 \lra v) ] \le \frac {C'' N^4} { |u-v|^{d-2}}
\end {equation}
(we write the right-hand sides in this way in order to emphasize that the interpretation is that when conditioning on the event $u \lra v$, the moments of the number of points on a self-avoiding connection from $u$ to $v$ correspond to a constant times $|u-v|^2$ points on this connection) that can be obtained easily from the BK inequality and the two-point function bounds. For instance, for the first sum, we can bound the sum separately depending on whether the distance of $t$ to $\{ u , v \}$ is smaller than $|u-v| / 3$ or not. We have clearly (summing over the value of the $L^\infty$ distance $k$ of $x$ to $u$)
\begin {eqnarray*}
 \sum_{t : \ |t - u| \le |u-v|/3}  P [ (u \circ t ) \circ (t \circ v) ] & \le & \sum_{k \le |u-v|<3 } \frac {C k^{d-1}}{k^{d-2}} \times \frac {1}{|u-v|^{d-2}}
 \le \frac {C' | u-v|^2}{|u-v|^{d-2}}
\end {eqnarray*}
(using the fact that $|v-t| >2 |v-u|/3$ in that case)
and the sum over $t$ that are not close to $u$ nor $v$ is easy to bound.
For the sum over $t_1$ and $t_2$, we can just use the previous bound for the sum over $t_1$ for each fixed $t_2$, and then (once the term $N^2$ has been factorized) use the same result again when summing over $t_2$. This type of sum appears for instance in \cite {CHS} where the authors estimate the moments of the chemical distances in critical percolation.

For each of the configurations (0), (1), (2')-(8'), when evaluating the probabilities fixing all the points except the two black dots (corresponding to $x'$ and $y'$) and summing over $x'$ and $y'$, we remove the two ``black points with just two incoming edges'' in the diagrams while adding a multiplicative factor $N^4$ in the probability of each one of the three diagrams \ref {g6}. For (9''), the same holds when summing over the two white dots. Then, (1) leads to (c), (2',5',6',7',8') to (a) and (0, 3',4',9'') to (b). The probability of the double connection (a) is clearly bounded by $|x-y|^{4-2d}$. But with $x$ and $y$ fixed at distance greater than $aN$ of each other, it is straightforward to check that for each of latter two diagrams (b) and (c), summing the probabilities over the inner points leads to an upper bound $O(N^{4 - 2d})$ as expected.

\begin{figure}[h]
  \centering
  \includegraphics[width=.6\textwidth]{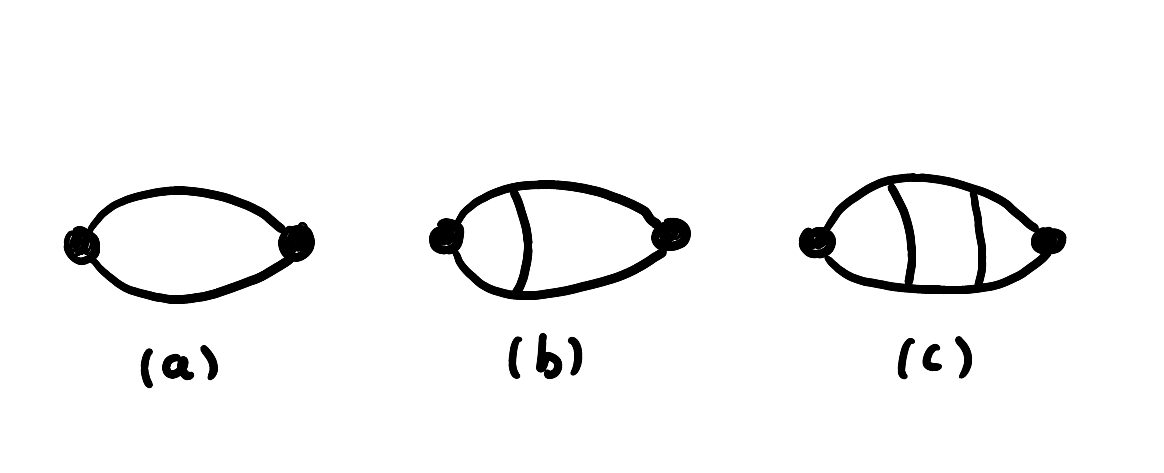}
  \caption{\label{g6}The diagrams after removing the middle points in edges}
  \end{figure}
\end {proof}

\subsection {Via lowest points}
\label {Sbrief}

Let now briefly outline how one can implement similar ideas to bound the second moment of the number $S$ of points in $\Lambda_N$ that have minimal first coordinate along an open self-avoiding loop of diameter greater than $aN$. The goal of this section is not to provide another detailed proof, but rather to illustrate how  the uniform half-plane connectivity bounds in \cite {panis_sharp_2025} allows to control a wide range of diagrammatic expansions.

\begin {itemize}
 \item As we have already explained, the uniform upper bound on $E[S]$ follows directly from Lemma \ref {L:restrict2}. In order to derive the lower bound, one can on the one hand rely on the estimates from \cite {panis_sharp_2025} that provide both a lower bound for $\Sigma_N$ (and in fact a variant of it in their Proposition 2.2, where one restricts the sum over point $t$ that are at distance $\eps N$ of the half-plane) and a uniform lower bound for the half-plane connection probability $P[ 0 \lra_H t]$. To conclude, one has to show that if one considers two independent samples $\omega$ and $\omega'$ for which the two events $0 \lra_{H \cap \Lambda_N} t$ and $t \lra_H 0$ are respectively satisfied (when $t \in \partial \Lambda_N$ is at $\eps N$-distance from the boundary of $H$), then the conditional probability that there exists an $\omega'$-open self-avoiding path from $t$ to $0$ in $H$ that stays at distance 1 from the cluster of the origin of $\omega$ (for percolation in $\Lambda_N$) if one removes $\Lambda_K$ from it is bounded from below. This can then be derived using the very same idea as in Lemma \ref {lemma1}.
 \item To upper bound the second moment, we need to sum over $x$ and $y$, the probability that $x$ and $y$ both contribute to $S$. The probability that they contribute disjointly is bounded via the BK inequality and leads to the upper bound $E[S]^2$ when summing over $x$ and $y$.  So, it remains to control the case where  $x$ and $y$ contribute without disjoint occurrence, i.e., when the two loops (that we will refer to as the $x$-loop and the $y$-loop) have to use a common edge. By symmetry, it clearly suffices to consider the case where the first coordinate of $x$ is not smaller than that of $y$.
  Following along the $y$-loop in both directions starting from $y$,  we can define the first points $u$ and $v$ at which the $y$-loop intersects the $x$-loop. Let us call these points $u$ and $v$.
 We are therefore in one of the two cases depicted in Figure \ref {fig:6b}. For fixed $x$, $y$, $u$ and $v$, the probability of such diagrams can be bounded by the BK inequality and the half-space connection estimates from \cite {panis_sharp_2025}. For the top connection between $u$ and $v$, one needs an additional estimate adapting the ideas leading to Lemma \ref {L:restrict2}.

\begin{figure}[h]
  \centering
  \includegraphics[width=
  \textwidth]{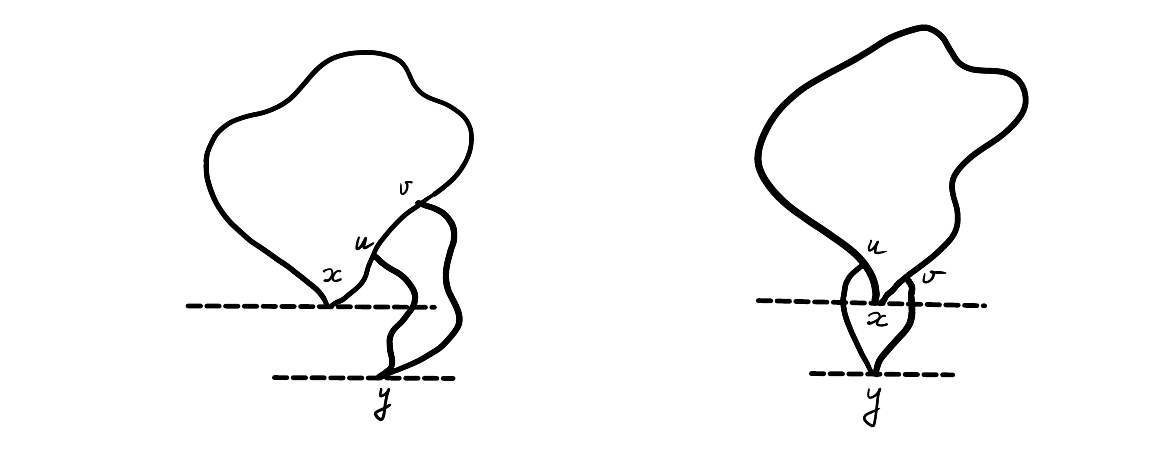}
  \caption{\label{fig:6b}The diagrams involved in bounding the second moment of $S$}
  \end{figure}

One can first estimate the case where $|u-x|$ or/and $|v-x|$ is/are greater than $aN/4$, which ends up providing a total contribution that tends to $0$ as $N \to \infty$. So, it remains to focus on the  case  where $|u-x|$ and $|v-x|$ are both smaller than $aN/2$.
Writing $\bar y = x-y$, $\bar u = u -x$ and $\bar v = v-x$, at the end of the day, one ends up for the contributions of the right-hand configurations in Figure \ref {fig:6b} with an upper bound looking like a constant times
$$ \sum_{\bar u, \bar v, \bar y \in (\N \times \Z^{d-1}) \setminus \{ 0 \} } \frac { 1}{ |\bar u|^{d-2} |\bar v|^{d-2} |\bar u+\bar y|^{d-1} |\bar v+ \bar y|^{d-1} }.$$
To check that this sum is finite, one can simply compare it with the corresponding integral in the continuum (and scaling immediately ensures that the integral indeed converges when $d >6$).
The left-hand configuration is treated similarly.
\end {itemize}

From there, one can also deduce Proposition \ref {mainprop} as in the previous section.

\section {Further results about the loops}

\subsection {$O(N^2)$ points per large loop-cluster}

We can combine/revisit the results and proofs in the previous sections to extract some information about the number of points on loop-clusters. More specifically:
\begin {enumerate}
 \item The expected number of points that belong to a loop-cluster of diameter greater than $aN$ in $\Lambda$ is bounded by a constant times $N^2$ (see the remark after Lemma \ref {L-horseshoe}).
 Hence, the probability that there exists a loop-cluster of diameter greater than $aN$ and with more than $bN^2$ points goes to $0$ as $b \to \infty$ uniformly in $N$.
 \item Lemma \ref {lemma3} states that $E[\pi^2]$ is bounded by a constant times $E[\pi]^2$, which implies by Paley-Zygmund (and the fact that $E[\pi]$ is bounded from below by a constant times $N^4$), that for some constant $c$, the probability
 that $\pi > c N^4$ is bounded from below by $c$.
 \item On the other hand, we have shown that the number of loop-clusters of diameter greater than $aN$ is tight. So, some fixed large enough $K$, the probability that there are more than $K$ different loop-clusters of diameter greater than $aN$ is smaller than $c/2$.
 \item Combining the last two items, we see that with probability at least $c/2$, there are no more than $K$ different loop-clusters of diameter greater than $aN$ and $\pi > cN^4$.
 So, at least one cluster must contribute  $cN^4 / K$ pairs of points to the sum $\pi$, which implies in particular that is has at least $N^2 (c/K)^{1/2}$ points and has diameter at least $aN$.
 It therefore follows that for sufficiently small but fixed $b=(c/ K)^{1/2}$, the probability that there exists (in $\Lambda_N$) a loop-cluster of diameter greater than $aN$ and with more than $bN^2$ points is bounded from below.
 \item Using the same argument (subdividing the box $\Lambda_N$ into finitely many smaller boxes) as in the final stretch of the proof of Proposition \ref {mainprop}, we can then conclude that for any $\eta$ and $k$, if one chooses $a$ and $b$ small enough,
 the probability that there exist at least $k$ distinct loop clusters in $\Lambda_N$, each with diameter greater than $aN$ and cardinality at least $bN^2$ is bounded from below by $1-\eta$.
\end {enumerate}

It is actually possible to derive a stronger result, namely that {\em any} large loop-cluster (i.e., of diameter at least $aN$) will have $O(N^2)$ points. More precisely:
\begin {proposition}
\label {P:length}
For any fixed $a$, and for any $\eta >0$, there exists $b_1$ and $b_2$ such that for all $N$, with probability at least $1- \eta$, any loop-cluster in $\Lambda_N$ with diameter at least $aN$ has at least $b_1N^2$ points and less than $b_2 N^2$ points
\end {proposition}

The upper bound (i.e., the statement with $b_2$) follows immediately from Point 1 above.
To derive the lower bound, one convenient option is to use the following relevant result adapted from the paper \cite {CHS} by Chatterjee, Hanson and Sosoe (this is statement (7) in their paper -- we could have also relied on earlier results of van der Hofstad and Sapozhnikov \cite {HStori}) -- we reformulate it so that it fits our setting: We let $E_{\lambda,N} (x,y)$  denote the event that there exists $\Lambda_N$ an open path joining $x$ and $y$ that has less than $\lambda |x-y|^2$ steps.
\begin {lemma}[\cite {CHS}]
\label {L:length}
There exists a constant $\xi$ such that for all $N$, for all $x$ and $y$ in $\Lambda_{N/4}$,
$$ P [ E_{\lambda, N} (x, y) ] \le \exp (-\xi / \lambda) \times P [ x \lra y].$$
\end {lemma}
The proof of the lemma goes via estimating the higher moments of the number of points lying on open paths from $x$ to $y$.

\begin {proof}[Outline of the proof of Proposition \ref {P:length} using Lemma \ref{L:length}]
For each point $x$ in $\Lambda_N$, we want to bound the probability of the event $U(x, N)$ that it is a lowest (i.e., the point with smallest first coordinate) point on a self-avoiding loop of $L^\infty$ diameter at least $aN$ with less than $\lambda N^2 / 64$ points. Note that such a loop will necessarily exit the box $\Lambda_{aN}$. We let $t$ denote the first point on $\partial \Lambda_{aN/4}$ visited by that loop and $t'$ denote the last point on $\partial \Lambda_{aN/8}$ visited by that loop (i.e., we can use the integer parts of $aN/4$ and $aN/8$ instead of $aN/2$ and $aN/4$).
\begin{figure}[h]
  \centering
  \includegraphics[width=.8\textwidth]{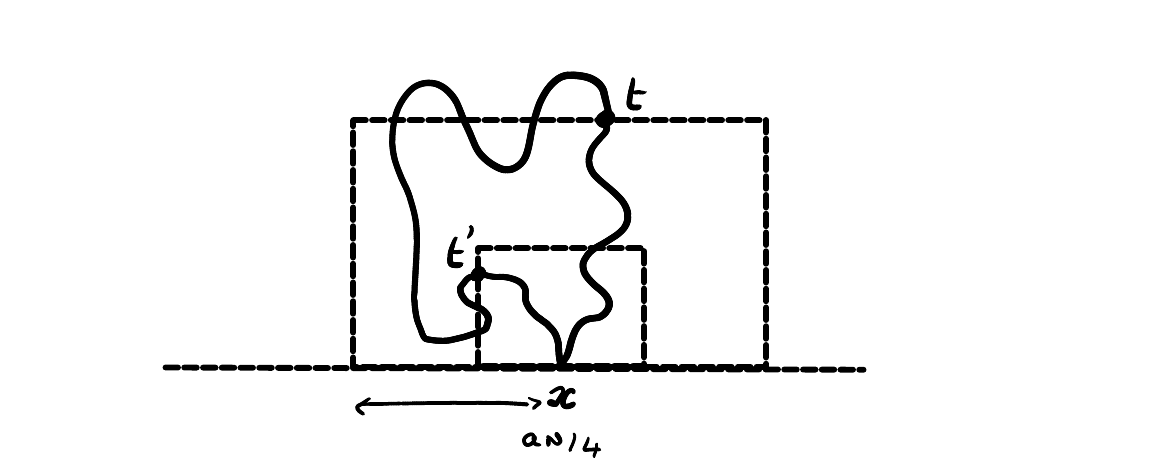}
  \caption{\label{g11}The decomposition of the loop -- the part between $t$ and $t'$ has to be a ``short connection'' between two points at distance at least $aN/8$ from each other}
  \end{figure}
The part of the loop between $t$ and $t'$ therefore joins two points that are at distance at least $aN/8$ from each other and it must have less than $\lambda N^2$ points. We therefore readily get an upper bound for $P[U(x,N)]$ of the type
$$ \Sigma_{aN/4}' \times \Sigma_{aN/8}' \times \sup_{t,t'} P [ E_{\lambda, N} (t, t') ] $$
which in turn is bounded by a constant (that depends on $a$ but not on $N$ or $\lambda$) times
$$\exp (- \xi / \lambda)  \times \Sigma_{aN/4}' \times \Sigma_{aN/8}' \times \sup_{t, t'} P [ t \lra t'] $$
which finally is bounded by a constant times
$$ \exp ( - \xi / \lambda) \times \frac {1}{N^d}.$$
Summing over all points $x$, we get that that the expectation of the number of such points is bounded by some constant times $\exp (- \xi / \lambda)$.
In particular, the probability that there exists one such point is bounded by this quantity independently of $N$. Since this quantity goes to $0$ as $\lambda \to 0$, this concludes the proof.
\end {proof}

\subsection {No pinching loops}

Let us now derive one further piece of information about these large loops:
We write $G(N,\eps, \alpha)$ for the event that there exist two points $x$ and $x'$ within distance smaller than $\eps N$ of each other such that in $\Lambda_N$, $(x\dlra^{\alpha N} x')\circ(x\dlra^{\alpha N} x')$ -- in this section and the next one,  $x \lra^k x'$ stands for $x \lra^{x + \Lambda_k} x'$).
This heuristically corresponds (for $\alpha$ fixed and $\eps$ very small) to the existence of a large loop that is $\eps N$ close to pinching.

\begin{lemma}\label{L:no-pinching}
    For fixed $\alpha<1$, the probability of $G(N, \eps, \alpha )$ is bounded by a constant times $\eps^{d-4}$ (uniformly in $N$).
\end {lemma}

\begin {remark}
This exponent $d-4$ can be interpreted in terms of analogous events for Brownian motions and Brownian loops in dimension $d$ (see the next section).
\end {remark}

\begin {proof}
A first remark is that  we can cover $\Lambda_N$ with $O(\eps^{-d})$ boxes of even side-length smaller than $\eps N$ in such a way that if $G(N, \eps, \alpha)$ occurs, then the  event $G(N, 4 \eps, \alpha/2 )$ occurs for some pair of points that are on the boundary of one of these boxes. When $B$ is such a box, we let $2B$ denote the box with twice its side-length and the same center point. By decomposing the loop into excursions from $\partial B$ to $\partial (2B)$, we see that if $G(N, \eps, a)$ holds, we can find one of these boxes $B$ and four points $y$, $y'$, $z$ and $z'$ on the boundary of $2B$ such that the following four events occur disjointly:
$$ (y \lra_{2B}^{c(B)} y'), \  (y' \lra^{\alpha N/2} z'), \ (z' \lra_{2B}^{c(B)} z), \ (z \lra^{\alpha N/2} y),$$
where $c(B)$ is the complement of $B$. In other words, the first and third events correspond to self-avoiding paths joining the boundary points of $2B$ within $2B$ that go all the way to the smaller box $B$.

\begin{figure}[h]
  \centering
  \includegraphics[width=\textwidth]{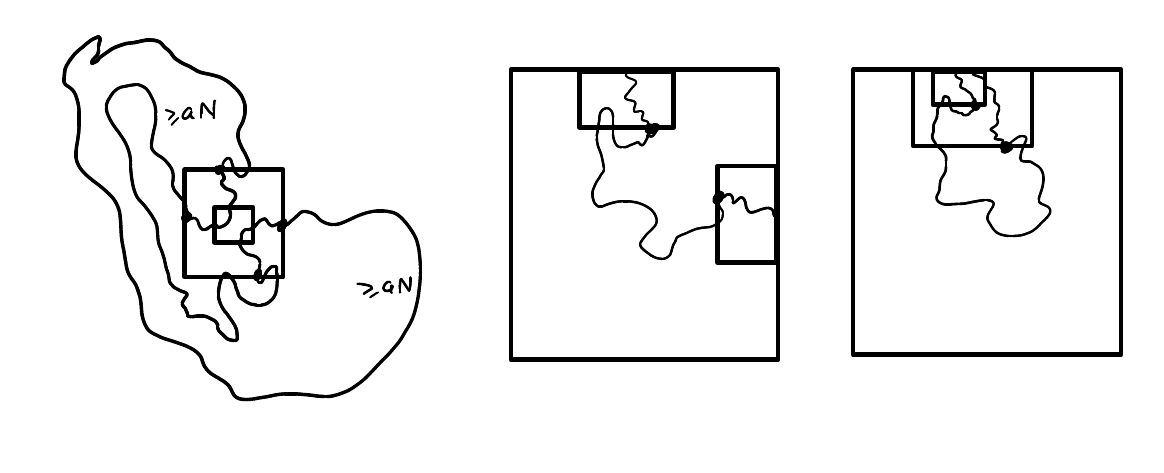}
  \caption{\label{g8}Decomposing the loop with the points $y,z,z',z$ (left), and the two ways to upper bound $P [ y \lra_{2B}^{c(B)} y']$ depending on whether $y$ are $y$ are not close to each other (left) or close to each other (right)}
  \end{figure}

The idea is that a self-avoiding loop will tend to go back and forth between $\partial B$ and $\partial (2 B)$ not more than a geometric number of times, so that the previous bound will not create too much overcounting.
It now suffices to use the BK inequality, to bound each term separately and to sum over all boxes and points $z, z', y, y'$.

The only missing ingredient is then provided by the following result that can be easily obtained using the same ideas as before, namely that there exists a constant $C$ such that
$$ P [ y \lra_{2B}^{c(B)} y'] \le C (\eps N)^{-d}.$$
Indeed:
\begin {itemize}
 \item If
$| y - y'| > \eps N / 100$, we can define $U$ (respectively $U'$) to be the intersection of $2B$ with the box of side-length $\eps N /400$ centered at $y$ (resp. centered at $y'$). We can then decompose the path from $y$ to $y'$ according to its first (resp. last) visited point on $\partial U$ (resp. $\partial U'$). Using Lemma \ref {L4}, we then readily obtain an upper bound of a constant times $((\eps N^{d-1})/(\eps N)^d ) \times (1/(\eps N)^{d-2}) \times ((\eps N)^{d-1} / (\eps N)^{d}) \asymp 1/(\eps N)^d $.
\item If $|y-y'| \le \eps N / 100$, we now choose ``concentric'' shells, i.e., we let $U$ (resp. $U'$) denote the intersection of $2B$ with the box of side-length $\eps N /10$ centered at $y$ (resp.
the box of side-length $\eps N /20$ centered at $y''$) and we decompose the path from $y$ to $y'$ (that goes through $B$, and therefore exits $U$ and $U'$) as before with the same bounds. This choice of shells ensures that the points on the different boundaries are at distance $\asymp \eps N$ from each other, so one can use the upper bound $1 / (\eps N)^{d-2}$ for the connection between them).
\end {itemize}
With this estimate in hand, we then see that for each of the $O(\eps^{-d})$ boxes $B$, the probability in question is bounded by a constant (that depends on $a$ but not on $\eps$ or $N$) times
$$ ((\eps N)^{d-1})^4 \times \frac {1}{N^{d-2}} \times \frac {1}{N^{d-2}} \times (\eps N)^{-d} \times (\eps N)^{-d}.$$
All the powers of $N$ cancel out and we indeed get the upper bound
$$ P [ G( N , \eps, \alpha) ] \le C \eps^{-d + 4d - 4 -d - d } = C \eps^{d-4}.$$
\end {proof}

\subsection {Loop-clusters are well-separated}

Almost the same proof as in the previous section allows to show that big loop-clusters remain at distance of order of $N^2$ apart.
More precisely, let us now consider the event
$G'(N,\eps, \alpha)$  that there exist two points $x$ and $x'$ within distance smaller than $\eps N$ of each other such that in $\Lambda_N$, $(x\dlra^{\alpha N} x)\circ( x' \dlra^{\alpha N} x' )$.
This corresponds (for $\alpha$ fixed and $\eps$ very small) to the existence of two large loop-cluster at distance less than $\eps N$ from each other (or the existence of two disjoint self-avoiding loops in the same cluster, which has a probability that goes to $0$ as $N \to \infty$ by Lemma \ref {L:no-glasses}).

\begin{lemma}\label{L:no-nearby}
    For all fixed positive $\alpha$, the probability of $G'(N, \eps, \alpha)$ is bounded by a constant times $\eps^{d-4}$ (uniformly in $N$).
\end {lemma}

\begin{figure}[h]
  \centering
  \includegraphics[width=\textwidth]{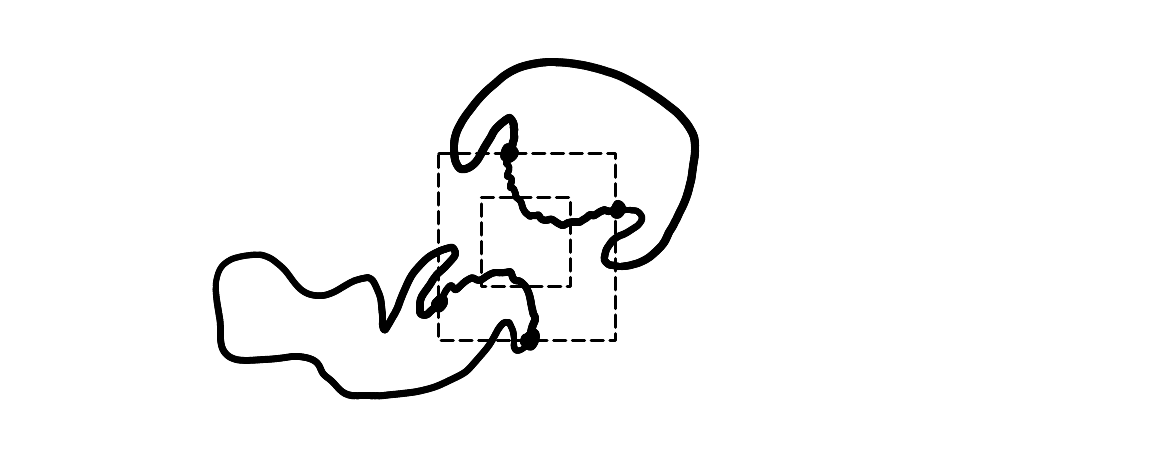}
  \caption{\label{z1}Decomposing the two loops with the points $y,z,z',z$}
  \end{figure}
The proof is almost identical to that of Lemma \ref {L:no-pinching}. The only difference is that in the decomposition with the excursions within $2B$, one now requires for the large connections of diameter at least $\alpha N/4$ that $z$ join $y'$ instead of $y$ and one requires $z'$ to join $y$ instead of $y'$ (so that one creates two loops instead of one) like in Figure \ref {z1} -- but all estimates and arguments remain unchanged, so we do not repeat it here.

\subsection {Scaling limits}

These last lemmas do in fact allow us to prove Proposition \ref {main3}. Let us first recall elementary facts about the topology of the space of self-avoiding loops in $[-1,1]^d$. We try to present this in an elementary way.

A self-avoiding loop $\gamma$ in $[-1,1]^d$ can be equipped with the Euclidean distance in $[-1,1]^d$. It is then a connected compact subset of $[-1,1]^d$, and it is also locally connected: For any $\alpha >0$, there exists $\delta >0$ such that for any two points $x$ and $y$ on $\gamma$ with $|x-y| < \delta$, the points are connected in the intersection of $\gamma$ with a ball of radius $\alpha$ [this follows immediately from the uniform continuity of $\gamma$].

A set $\Gamma$ of self-avoiding loops in $[-1,1]^d$ is said to be uniformly locally connected if for any $\alpha$, one can choose a $\delta$ that would fit all $\gamma$'s in $\Gamma$ simultaneously.
It turns out that this condition is essentially ensuring tightness in the set of self-avoiding loops: Suppose that one has a sequence $(\gamma_i, i \ge 1)$ of self-avoiding loops of diameter at least $a$ in $[-1,1]^d$ that is uniformly locally connected, then one can find a subsequence $i_n \to \infty$ such that $\gamma_{i_n}$ converges (for the Hausdorff topology on compact sets) to a self-avoiding loop $\gamma$. This is a rather classical result: First one can use the fact that the set of all compact subsets of $[-1,1]^d$ is compact for the Hausdorff topology, so that there exists a subsequence that converges to some compact set $\gamma$ of diameter at least $a$; it remains to check that $\gamma$ is a self-avoiding loop. For this, we note that $\gamma$ is necessarily connected, and that the uniform local connectedness assumption immediately implies that one can use the same $\delta=\delta(\alpha)$ to ensure that $\gamma$ is also locally connected.
Now we want to use Whyburn's Theorem \cite {Wh,Nadler} that states that a continuum set $\gamma$ is a self-avoiding loop (i.e., homeomorphic to the unit circle) if and only if it has no cut points and if removing any pair of points on $\gamma$ disconnects it. Both conditions are easy to check here: If three points $x$, $y$ and $z$ are in $\gamma$, then one can find $x_{i_n}$, $y_{i_n}$ and $z_{i_n}$ in $\gamma_{i_n}$ that converge to them as $n \to \infty$, and then, by looking at the (subsequential) limit of the part of $\gamma_{i_n}$ that does not contain $x_{i_n}$, we see that $y$ and $z$ are still connected within $\gamma \setminus \{ x \}$, so that $x$ is not a cut point of $\gamma$. Conversely, for each $\eps$, by observing that the intersections of the two parts of $\gamma_{i_n} \setminus \{ x, y\}$ with the complement of the $\eps$-balls around $x$ and $y$ remain separated, we see that the same is true for the two parts of $\gamma$. Since this is true for all $\eps$, one can deduce that $\gamma \setminus \{ x, y \}$ is not connected.

\begin {proof}[Proof of the proposition]
Since the space of compact subsets of $[-1,1]^d$ is compact for the Hausdorff topology, it immediately follows that  for any sequence $N_n' \to \infty$, there exists a subsequence $N_n \to \infty$ so that for all $j \ge 1$, the joint law of $(L_1(N_n), \ldots, L_j(N_n))$ converges (weakly with respect to the Hausdorff topology). This induces a limiting law on infinite sequences of compact subsets of $[-1, 1]^d$.  Let us denote by $(L_1, \ldots, L_j, \ldots)$ a sample of the limiting law. It follows readily from Proposition \ref {mainprop} that if $\rho_j$ denotes the diameter of $L_j$, then $\rho_j$ is almost surely a sequence of positive numbers that tends to $0$ as $j \to \infty$. Furthermore, Lemma \ref {L:no-nearby} readily implies that any two $L_j$'s are almost surely disjoint. Indeed, for all fixed $a$, the probability that two $L_j$'s of diameter strictly greater than $a$ in this family are less than $\eps$-apart goes to $0$ as $\eps \to 0$ (this follows by passing the statement Lemma \ref {L:no-nearby} to the limit), so that any two $L_j$'s  of diameter greater than $a$ are almost surely disjoint (and we can then let $a \to 0$). It therefore only remains to show that each $L_j$ is almost surely a self-avoiding loop, which is where Lemma \ref {L:no-pinching} enters the game.

Let us fix $j$ and $\eta >0$.
Lemma \ref {L:no-pinching} shows that for  any $\alpha_k = 2^{-k}$, there exists $\eps_k$, such that for any $n$, with probability at least $1 - \eta 2^{-k}$ the following is true: Any two points on $L_j(N_n)$ that are less than $\eps_k$ apart are joined by an arc of $L_j(N_n)$ of diameter smaller than $\alpha_k$. Hence, for any $n$, with probability $1-\eta$, this property holds simultaneously for all $k$. This in particular shows that there exist a set of uniformly locally connected self-avoiding loops $S$ (corresponding to these choices of $\eps_k$'s) such that for each $n$, the probability that  $L_j(N_n)$  belongs to $S$ is at least $1- \eta$. This implies immediately that with probability at least $1- \eta$, the same holds for $L_j$ --- which by the previous analytical topology consideration implies that $L_j$ is a self-avoiding loop with probability at least $1- \eta$. Since this is true for all $\eta$, $L_j$ is almost surely a self-avoiding loop.
\end {proof}

\section {Other events with scaling limits}
\label{S73}

\subsection {Starfishes}
\begin{figure}[h]
  \centering
  \includegraphics[width=
  \textwidth]{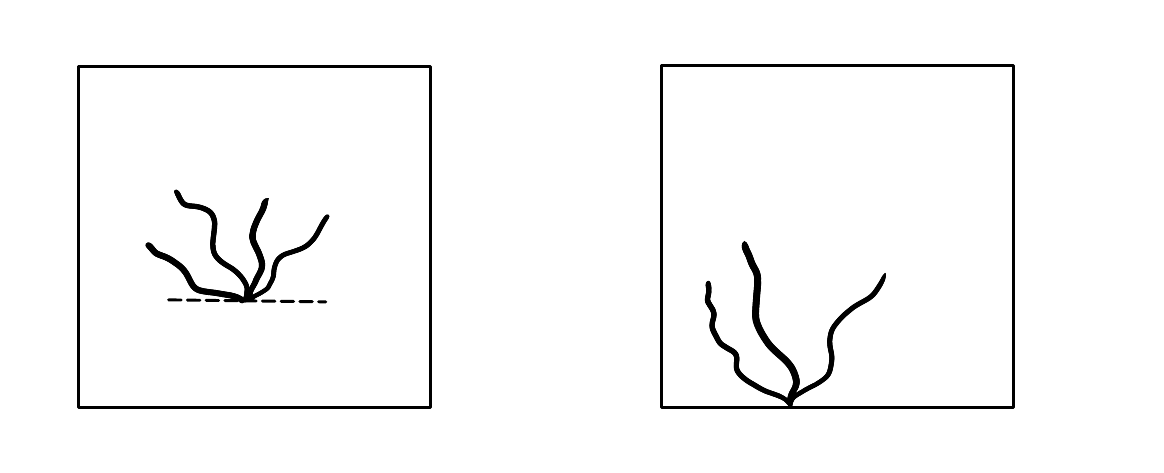}
  \caption{\label{p2}Restricted starfishes}
  \end{figure}

We now describe some types of large clusters that are tight and have scaling limits (here we always assume that $d$ is large enough so that the two-point estimate (\ref {tpe}) holds):

\begin {itemize}
\item
A first natural class comes when one looks for ``star-shaped'' configurations within clusters:
 One can look at starfishes with $l$ arms for $l \ge 4$. We can for instance fix $a<1$ and look at the set of points $z$ such that $z \lra z+ \partial \Lambda_{aN}$ occurs disjointly $l$ times. This corresponds to the case depicted  in Figure \ref {p1}. The expected number $n(a,N)$ of such points in $\Lambda_N$ can be shown to be  comparable to $N^d N^{-2l}$  (here, the lower bound comes from the one-arm estimates from \cite {KN} and one can use the same argument described in Section \ref {Sdc} to obtain the lower bound). This quantity will therefore be $\asymp 1$ when $d=2l$.
\item
Variants include starfishes with additional restrictions:
For instance, one can look at the set of points $ z $ on the boundary of $\Lambda_N$ for which one can find the event $z \lra_{\Lambda_N} z + \partial \Lambda_{aN}$ occurs disjointly $l$ times (i.e., one can think of an $l$-arm starfish or plant attached at the side of the $d$-dimensional aquarium). Then one has $N^{d-1}$ a priori options for $z$ and each arm costs a factor $N^3$ (see the half-space exponent in \cite {CH}), so that the number of such $z$ is tight when $d=3l + 1$. The right sketch in Figure \ref {p2} is the case $l=3$.
\item One can also consider half-space starfishes inside of $\Lambda_N$, i.e., look for the set of points $z \in \Lambda_N$ for which $z \lra_{z+H} \partial \Lambda_{aN}$ occurs disjointly $l$ times, where $H$ is a half-space going through the origin. This time, one has $N^d$ options for $z$, so that the number of such points $z$ will be tight when $d=3l$.
\item
One can look for other type of restrictions -- for instance  cone points for starfishes etc.
\end {itemize}

In all these examples, the fact that the expected number of such points $z$ in $\Lambda_N$ is bounded comes from the existing estimates for one-arm probabilities. One strategy to show the existence of such configurations is  to derive a matching lower bound on this expectation and to also get an upper bound on the second moment of such points $z$.

One important remark is that as opposed to the loop-case, where the loop-cluster allowed to identify essentially one loop at macroscopic level, the starfishes within a cluster are not unique. Indeed, the shape of such starfishes is more like that of $l$ disjoint trees (and one can pick one out of several possible branches in each tree to draw the starfish) that are joined at their root (see Figure \ref {p1}).

Let us outline in some detail how the proof goes in the case of starfishes in the bulk. Let us start with the proper definitions:

\medbreak \noindent
{\sc Definition.}
{\em
We say that a percolation cluster contains an $l$-arm starfish of radius $N$ if it contains a point $z$ for which $z \lra z + \partial \Lambda_{N}$ occurs disjointly $l$ times.
We say that a point is a center of an $l$-arm starfish of radius $N$ if $z \lra z +  \partial \Lambda_{N}$ occurs disjointly $l$ times. So, an $l$-arm starfish contains at least one center.}
\medbreak

We now consider percolation in the whole of $\Z^d$, choose $a>0$ and define $n(a,N)$ the number of $l$-arm starfishes of radius $aN$ with a center in $\Lambda_N$.
\begin {proposition}
\label {starfish}
Let $l \ge 4$ and $d = 2l$. If the two-point estimate (\ref {tpe}) holds in dimension $d$, then:
\begin {itemize}
\item
The laws of the random numbers $n(a,N)$ are tight.
\item And on the other hand for any $k \ge 0$ and $\eta >0$, there exists $a>0$ such that $P[ n(a,N) \ge k ] \ge 1- \eta$.
\end {itemize}
\end {proposition}

The first step is to control the first moment of the number $m (a,N)$ of centers of starfishes of radius $aN$ that lie in $\Lambda_N$. We work always under the assumption that $l \ge 4$ and that (\ref {tpe}) holds for $d=2l$.
\begin {lemma}
For each fixed $a>0$, one has $E[ m (a, N)^2] \asymp 1$.
\end {lemma}
\begin {proof}[Outline of the proof]
Recall the one-arm estimate from \cite {KN} (see also \cite {panisschapira2})
$$ P [ 0 \lra \partial \Lambda_n ] \asymp n^{-2}.$$
The BK inequality (and summing over the points $z$ in $\Lambda_N$ immediately implies the upper bound for $E[ m(a,N)]$. To derive the lower bound, one needs to check that
$$
w_n := P [ (0 \lra \partial \Lambda_n) \hbox { occurs disjointly }   l   \hbox { times}] \ge c P [ 0 \lra \partial \Lambda_n]^l.$$
This can be achieved using the exact same strategy as in Lemma \ref {lemma1}, by comparing the picture with that of $l$ independent percolation samples.
\end {proof}
This lemma already shows that $m(a, N)$ is tight.
The next step is to bound the second moment of $m(a,N)$, i.e., to show that
\begin {lemma}
One has $ \sup_N E[ m(a,N)^2 ] < \infty$.
\end {lemma}
\begin {proof}[Outline of the proof]
We are going to use a diagrammatic approach. Suppose that $y$ and $z$ are two points in $\Lambda_N$.
Let us denote by $L(y,z)$ the event that  they are both centers of $l$-arm starfishes of radius $aN$.
We can decompose the event $L(y,z)$ according to the maximal number $\lambda$ for which the event that $z$ is the center of an $l$-arm starfish of radius $aN/2$ and the event that $y$ is the center of an $\lambda$-arm starfish of radius $aN/2$ occur disjointly. In other words, if $\lambda=l-k$, then at least $k$ arms of the starfish centered at $y$ have to intersect a starfish centered at $z$.

The BK inequality immediately shows that
$$ P [ L(y,z), \lambda= 0 ] \le w_{aN/2}^2 \le C N^{-2d}$$
for some constant $C$. Summing over all $y$ and $z$ in the box $\Lambda_N$ gives a constant upper bound.

\begin{figure}[h]
  \centering
  \includegraphics[width=.8
  \textwidth]{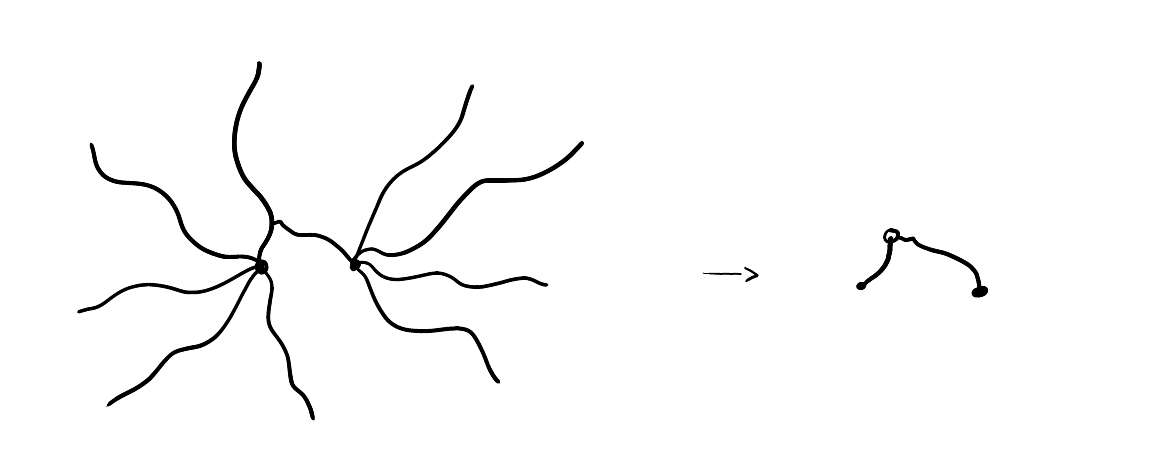}
  \caption{\label{h1}Two joined $l$-starfishes with a total of $2l-1$ disjoint arms, and the corresponding trimmed subdiagram}
  \end{figure}
When $\lambda = 1$, one has a diagram as on the left of Figure \ref {h1}, and using the BK inequality and the one-arm probability upper bound for each of the $2l-1$ arms of radius at least $aN/2$, we get the upper bound
$$ P [ L(y,z),\  \lambda= 1] \le C N^{-2(2l-1)} \sum_t P [ (z \lra t) \circ (t \lra y) ] \le C' N^{-4l + 2} \frac {1}{|z-y|^{d-4}},$$
where we used (\ref {final1}) in the latter step to sum over the position of the white dot in the right of Figure \ref {h1}. Summing over all $z$ and $y$ in $\Lambda_N$ then gives an upper bound $C / N^{d-6}$.

\begin{figure}[h]
  \centering
  \includegraphics[width=.8
  \textwidth]{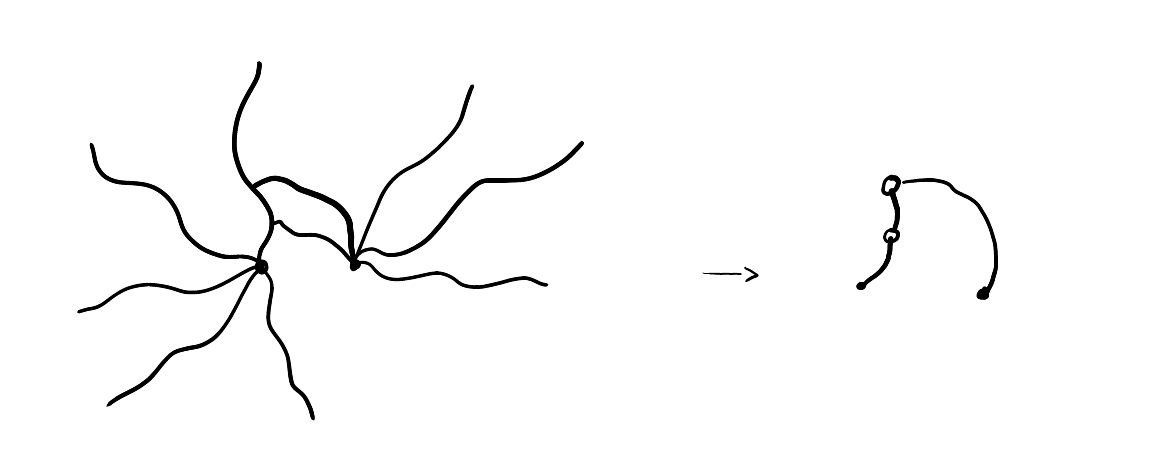}
  \caption{\label{h2}Two doubly-joined $l$-starfishes with a total of $2l-2$ disjoint arms (option 1) and the corresponding trimmed subdiagram}
  \end{figure}

\begin{figure}[h]
  \centering
  \includegraphics[width=
  .8\textwidth]{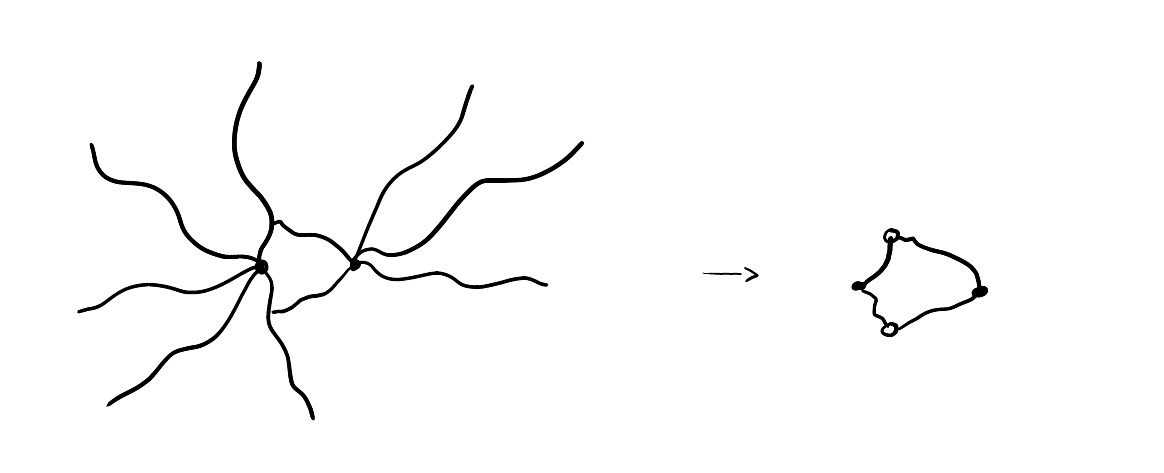}
  \caption{\label{h3}Two doubly-joined $l$-starfishes with a total of $2l-2$ disjoint arms (option 2) and the corresponding trimmed subdiagram}
  \end{figure}
When $\lambda=2$, depending on whether the two arms starting at $y$ that touch the arms starting at $z$ on the same arm or not, one ends up with one of the configurations depicted on Figure \ref {h2} or Figure \ref {h3}. We can then use the same procedure as before by first bounding the probability of the $2l-2$ arms in the diagram by a constant times $N^{-2(2l-2)}$, and it suffices to bound the two sums over the trimmed (i.e., one has cut off the arms) diagrams depicted on the right-hand side of these two figures (note that for the first one, we simply dropped a connection as it does converge without it anyway):
$$ \sum_{z,t_1, t_2, y} P [ (z \lra t_1) \circ (t_1 \lra t_2) \circ (t_2 \lra y) ] \hbox { and }
\sum_{z, t_1, t_2, y} P [ (z \lra t_1) \circ (t_1 \lra y) \circ (z \lra t_2) \circ (t_2 \lra y) ] .
$$
Both can be bounded as above by using (\ref {final2}) and (\ref {final1}) and lead to an $o(1)$ bound.

\begin{figure}[h]
  \centering
  \includegraphics[width=.8
  \textwidth]{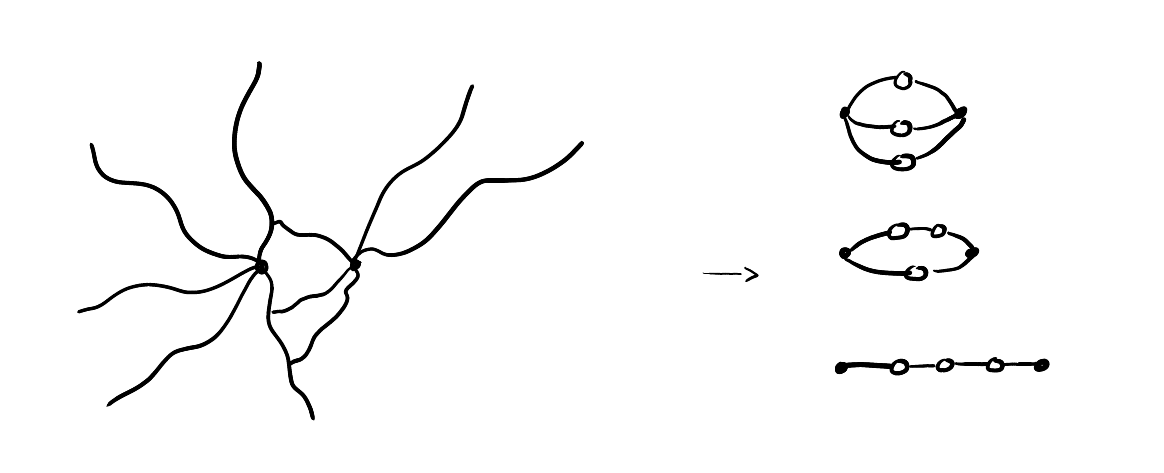}
  \caption{\label{h4}Two triply-joined $l$-starfishes and the trimmed subdiagrams}
  \end{figure}
For the final case $\lambda \ge 3$, we can bound things brutally by only focusing on three arms from $y$ that touch the $l$ arms from $z$, and forget about the other $l-3$ ones (indeed, in this way, $y$ will be a trifurcation for this diagram, which will ensure the convergence of the sum).
Then, depending on if these three arms land on similar arms or not, one ends up summing over the three diagrams depicted on the right of Figure \ref {h4} (the one corresponding the left-hand side of the figure would be the middle diagram on the right). All three sums are easily shown to converge (the third one requires the generalization of (\ref {final2}) to three intermediate point), which when summing over $z$ (which compensates the bounds on the $l$ arms) one readily checks that
$$ \sum_{z,y} P [  L(y,z),  \ \lambda\ge 3 ] \le o(1)$$
because $d \ge 7$.

So, wrapping up, we get indeed a constant upper bound for $E[ m(a,N)^2]$.
\end {proof}

\begin {proof}[Outline of the proof of the Proposition]
We can then conclude as in the case of the loop-clusters:
 The two previous lemmas and Cauchy-Schwarz show that the probability that $P[ m(a,N) \ge 1]$ is bounded from below by a positive constant $c(a)$ when $a$ is small enough. We fix such an $a<1/2$.
 We then choose $l$ so that the sum of $2^{ld}$ independent Bernoulli random variables of parameter $c(a)$ is greater than $k$ with probability at least $1- \eta/2$. We then choose $N_1 > N 2^{-l} /4$ such that one can  fit $2^{ld}$ disjoint boxes of the same size as $\Lambda_{2N_1}$ in $\Lambda_N$. In each of these boxes, the probability that there exists an $l$-arm starfish with radius $aN_1$ in its middle box (the box with same center, but scaled by a factor $1/2$) is at least $c(a)$. The events in these $2^{ld}$ boxes are independent, so we conclude that with a probability at least $1 - \eta/2$, at least $k$ of these will contain such an $l$-arm starfish.

 It remains to check that these $l$-arm starfishes are necessarily in different clusters. This is quite easy to do because the starfishes necessarily lie in different boxes and are  at least $N_1/2$-distant from each other, so that the expected number of pairs of starfish centers in the different boxes will be bounded by a constant times $(N^2 N^2) /N^{d-2} = 1/N^{d-6}$ and therefore goes to $0$.
\end {proof}

\subsection {Further  shape types}
\label {Sfurthershapes}

\begin{figure}[h]
  \centering
  \includegraphics[width=
  \textwidth]{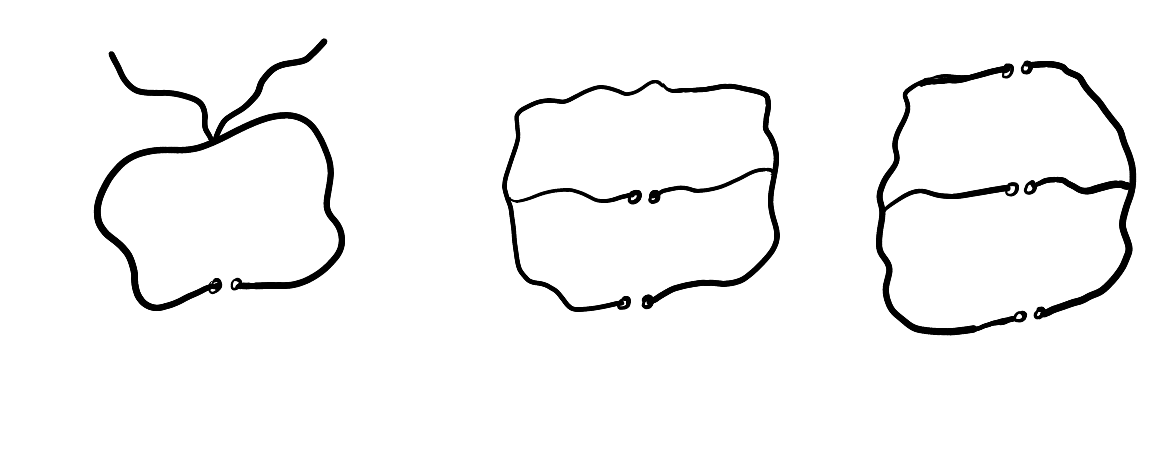}
  \caption{\label{p3} Some configurations of interest (in any dimensions for the left-hand type, in $d=10$ for the middle one, in $d=12$ for the one on the right)}
  \end{figure}

Another class of rather natural shapes to look at involve more than one special marked point (in the starfishes, one could look for its center, in the cases of loops, one could look for the lowest point).
Among graphs with two marked points, one can for instance consider the case depicted on the left of Figure \ref {p3}. We can fix $a<1$ and look for the set of pairs of points $x$ and $y$ at distance at least $aN$ from each other such that the following four events occur disjointly:
 $x \lra y$, $x \lra y'$ where $y'$ is a neighbor of $y$, $x \lra x+ \partial \Lambda_{aN}$, $x \lra x+ \partial \Lambda_{aN}$.
 Again, by combining the same ideas, one easily sees that if $n(a,N)$ denotes the number of such pairs,
 $$ E[ n(a,N)] \asymp N^{2d} \times \frac 1 { N^{d-2} N^{d-2} N^2 N^2 } \asymp 1$$
 regardless of the dimension.

 We can also look at more complicated situations with more than one cut.
For instance, one can consider the set of quadruples of points $(z_1, z_2, x_1, y_1)$
at distance at least $aN$ from each other for which for some neighbour $x_2$ from $x_1$ and some neighbour $y_2$ from $y_1$ (see the middle part of Figure \ref {p3}), one has
$$ ( z_1 \lra z_2 ) \circ ( z_1 \lra x_1 ) \circ ( z_1 \lra y_1 ) \circ ( z_2 \lra x_2 ) \circ ( z_2 \lra y_2 ).$$
The expected number of such collection of points $(z_1, z_2, x_1, x_2)$ can then be seen to be
$$ \asymp N^{4d} \frac {1}{N^{5(d-2)}} \asymp N^{10-d}.$$
So, we see that the case $d=10$ is special here, provided the two-point function estimate is also valid in this case.

One can also look at configurations cxreated by pairs of neighboring clusters instead of single clusters. For instance, one can look at pairs of three-arm stars that are ``glued'' pairwise by each of their three arms as in the right sketch of Figure \ref{p3}. This corresponds to the case of points $(z_1, z_2, x_1, y_1, t_1)$ in $\Lambda_N$ that are at distance greater than $aN$ from each other, so that for some neighbours $x_2, y_2, t_2$ of $x_1, y_1, t_1$ respectively,
$$ (z_j \lra x_j) \circ (z_j \lra y_j) \circ (z_j \lra t_j)$$ holds disjointly for both $j=1$ and $j=2$. One has five marked points and six connections, so the expected number of such quintuples will be $\asymp N^{5d}/N^{6(d-2)} \asymp N^{12-d}$. So, this type of configuration will be of interest when $d=12$.
\medbreak

Again, one can also add special restriction  (such as requiring the cluster to have a global cone point) or to ask for special points to be located on the boundary of the box.
\medbreak

In each case, the strategy to derive the tightness will be the same, but the second moment bounds will be a little bit more involved, due to the presence of more than one marked points, which leads to more diagrams. One also has to be wary that when looking at neighboring clusters or arms  then one will need a $d \ge 9$ (or at least $ \ge 8$) condition to ensure that the typical number of pairs of touching points remains tight (so, this restriction on $d$ is  needed for the first type depicted in Figure \ref {p3}).

\section {Brownian limits}

We make some final remarks that we hope to expand upon in  \cite {CW}.

\subsection {Towards loop-soups}

As we have already mentioned, the probability of occurrence of individual loop-clusters should converge to a multiple of the Brownian loop-measure in the scaling limit, in a similar way as the IIC backbone  converges to Brownian motion, see \cite {HS2,HS3,KN2,HHHM} for related estimates and results. This is also related to the recent work of Chatterjee, Chinmay, Hanson and Sosoe \cite {CCHS} about the distribution of the chemical distance between connected points.

In addition to this, when one conditions on the existence of a given loop-cluster in some small tube, then the rest of the configuration outside of the tubes is not affected. This (together with estimates that show that any two big loop-clusters will be at macroscopic distance from each other) will then make it possible to deduce that any (sub)-sequential limit of the collection of all big loop-clusters would necessary be a Poisson point process, i.e., a Brownian loop-soup \cite {LW} if one has good control of the intensity measure.

As we have already mentioned, the exponent $\eps^{d-4}$ at the end of the previous section would then be interpreted as the  probability (for fixed $a$ and very small $\eps$ that a $d$-dimensional continuous Brownian loop  $\gamma$ (of time-length $1$ say) has a ``close-to-double-point'', i.e. that for four times $s_1, s_2, s_3, s_4$ circularly ordered, one has $|\gamma(s_2) - \gamma (s_4)| < \eps$ while $\gamma (s_1)$ and $\gamma (s_3)$ are both at distance greater than $a$ from $\gamma(s_2)$.

Similarly, when the dimension of space is large enough, the scaling limit of all the other tight configurations described in the previous section should be described in terms of Poisson point processes of  corresponding ``Brownian figures'' (the starfishes with $l$ arms would for instance correspond to $l$ independent integrated super-Brownian excursions that are glued together at the center).

\subsection {Almost-loops, almost-almost-loops}
\begin{figure}[h]
  \centering
  \includegraphics[width=
  \textwidth]{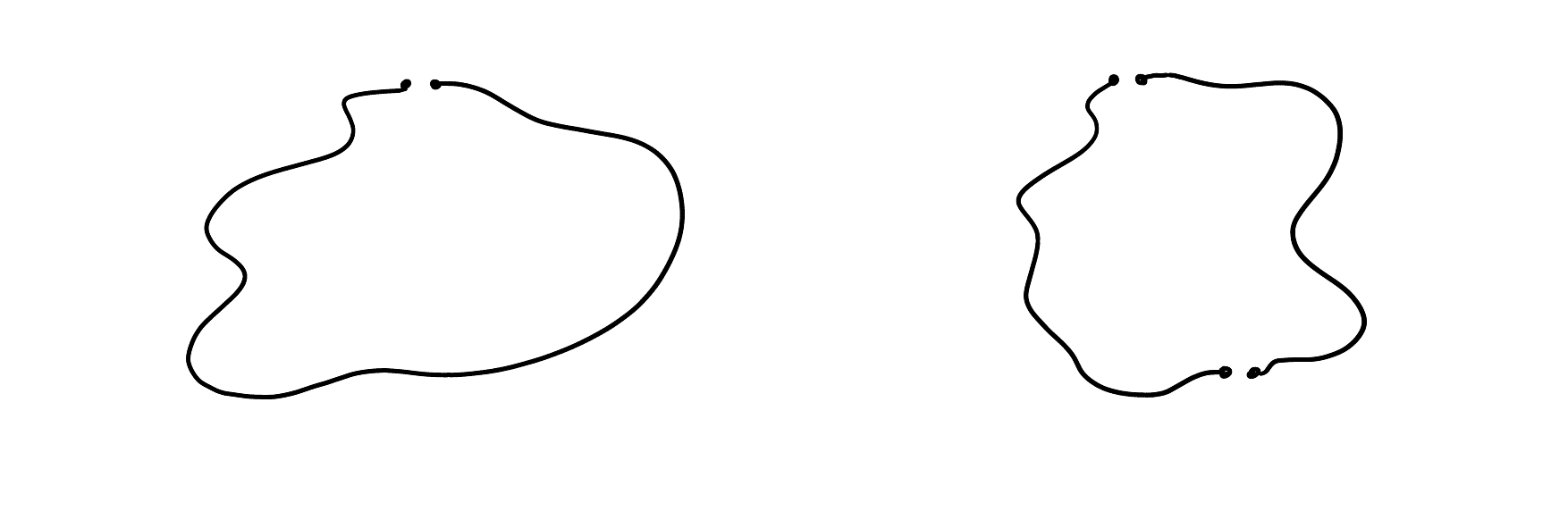}
  \caption{\label{fig:7} An almost-loop and an almost-almost-loop.}
  \end{figure}

Let us consider dynamical percolation where the status of different edges are updated independently after exponential waiting times. A loop cluster will have many ``cut-edges'' (such that the removal of the edge will disconnect the loop). Conversely, when one has an ``almost-loop-cluster'' (a cluster with only one missing edge to contain a large loop, as on the left of Figure \ref {fig:7}), then it will typically have only very few possible edges that would create a big loop (provided $d \ge 9$). So, the probability that a loop-cluster gets dislocated is much larger than that of a given almost-loop cluster turning into a loop-cluster. This shows that the number of almost loop-clusters will have to be much larger than then number of actual loop-clusters. In fact, since a typical large loop in $\Lambda_N$ would have of order $N^2$ points, one can infer  that the number of such almost loop-clusters will be of the order of $N^2$ when $d >8$. So, while these almost loop-clusters are very exceptional among the circa $N^{d-6}$ large clusters, they will nevertheless proliferate. The almost loop-clusters will also be Brownian-like (i.e. the clusters containing loop-clusters but with one marked disconnection point on the loop) in the scaling limit.
Similarly, we can look for configuration of loop-clusters with two pivotal edges erased, i.e. ``pairs of large clusters'' that do neighbor each other at two far-away points (so that the union of the two clusters with two additional edges would create a large loop, as on the right sketch in Figure \ref {fig:7}). The number of those pairs will then be of the order of $N^4$ (again, provided $d$ is not too small), and we can view these pairs as Brownian loops with two marked disconnection points on them.

\subsection {Near-criticality and massive loop-soups}
We now return briefly to the question of the scaling limit in the ``near-critical window''. Near-critical percolation in high-dimensions has been well-studied and is now well-understood (see for instance \cite {HH,CHS,HMS} and the references therein). In our setting, if one takes $p=p_c - u N^{-2}$ and considers critical percolation within the box $\Lambda_N$, a large loop-cluster for $p_c$ will anyway survive at $p$ (for the natural coupling between the two percolation realizations) as long as none of the $O(N^2)$ edges forming one of the big loops has been closed. On the other hand, the loop-cluster will disappear (i.e., not be macroscopic anymore) if any of the $O(N^2)$ ``pivotal'' edges has been flipped. This leads naturally to the idea (that we again plan on detailing in future work) that the scaling limit of collection of loop-clusters will then be a massive Brownian loop-soup (i.e., one starts with the usual Brownian loop-soup in $[-1,1]^d$ and independently erases each Brownian loop with a probability $e^{-m T(\gamma)}$ for $m = cu$ for some constant $c$, where $T(\gamma)$ is the time-duration of the loop). A similar feature would  work also for positive $u$.

\subsection*{Acknowledgments}
This research has been supported by grants from the Royal Society, the Cambridge Trust and a NSERC scholarship. We thank an anonymous referee for insightful comments on an earlier much shorter (and very different) version of this paper.

\printbibliography

\end{document}